\documentclass{amsart}
\usepackage{amssymb,amsmath,amsfonts,amsthm}
\usepackage{graphicx}

\input diagrams

\newtheorem{thm}{Theorem}[section]
\newtheorem{prop}[thm]{Proposition}

\newtheorem{cor}[thm]{Corollary}
\newtheorem{lemma}[thm]{Lemma}
\newtheorem{defn}[thm]{Definition}
\newtheorem{conj}[thm]{Conjecture}

\theoremstyle{definition}
\newtheorem{ex}[thm]{Example}
\newtheorem{rem}[thm]{Remark}

\errorcontextlines=0 \numberwithin{equation}{section}
\def\A{\mathcal{A}}
\def\a{\alpha}

\def\b{\beta}
\def\C{\mathbb{C}}
\def\ca{\check{\a}}
\def\cb{\check{\b}}
\def\cent{\mathcal{Z}}

\def\e{\epsilon}
\def\F{\mathbb F}
\def\g2#1#2{\setlength{\unitlength}{1cm}
\begin{picture}(3,1)
\put(1,0.1){\circle{0.2}} \put(1,0){\line(1,0){1.2}}
\put(1.1,0.1){\line(1,0){1}} \put(1,0.2){\line(1,0){1.2}}
\put(2.2,0.1){\circle{0.2}} \put(1.6,0){$\rangle$}
\put(0.9,0.3){#1} \put(2.1,0.3){#2}
\end{picture} }
\def\G{\mathcal{G}}
\def\ga{\gamma}
\def\H{\mathcal{H}}
\def\half{\frac{1}{2}}
\def\Hrcc{\hat{\H}^t_\R}
\def\Hgr{\mathbb{H}}
\def\Hgrrcc{\hat{\Hgr}^t_\R}
\def\I{\mathcal{I}}

\def\l{\lambda}

\def\O{\mathcal{O}}

\def\P{\mathcal{P}}

\def\pf{\noindent{\it Proof}}
\def\qed{\hfill $\square$\par\medskip}
\def\Q{\mathbb{Q}}

\def\R{\mathbb{R}}
\def\real{\mathfrak a}

\def\S{\mathcal{S}}

\def\t{\mathfrak{t}}
\def\trm{\mbox{tr}_m-\mbox{Ind}}

\def\Z{\mathbb{Z}}

\begin{document}

\begin{abstract} Recently Delorme and Opdam have generalized the theory of $R$-groups towards affine Hecke algebras with unequal labels.
We apply their results in the case where the affine Hecke algebra is
of type $B$, for an induced discrete series representation with real
central character. We calculate the $R$-group of such an induced
representation, and show that it decomposes multiplicity free into
 $2^d$ irreducible summands. The power $d$ can be calculated combinatorially.
\end{abstract}

%\begin{abstract} Let $\H$ be an affine Hecke algebra on a root system of type $B_n$ with possibly unequal parameters. In \cite{S}, we have formulated a conjecture which describes a parametrization of the set $\Hrcc$ of equivalence classes of irreducible tempered representations of $\H$ which have real central character. The non-discrete series representations in $\Hrcc$ are all obtained as irreducible components of the induction of a discrete series representation of a Levi algebra associated to a standard parabolic root subsystem. The conjecture also predicts the number of irreducible components of such an induced representation. Recently, Delorme and Opdam (\cite{DO2}) have given a description of a group (``the $R$-group'') which determines the structure of the intertwining algebra of these induced representations. We use this general description to compute explicitly and uniform in the parameters the $R$-group for $\Hrcc$ as above (in terms of certain Young tableaux). This description yields the number of irreducible components of the above induced representations. We obtain a multiplicity one decomposition. This result agrees with the above mentioned conjectures.\end{abstract}
\author{K. Slooten}
\email{slooten@ihes.fr}
\date{\today}
%\keywords{brol}
%\subjclass[2000]{primary: ...; secondary: ...}
\address{European Post-Doctoral Institute\\IHES\\Le Bois-Marie, 35, Route de Chartres\\ 91440 Bures-sur-Yvette\\France}
\title[The $R$-group for type $B$ Hecke algebras]{Reducibility of induced discrete series representations for affine Hecke algebras of type $B$}
\thanks{During part of the research for this paper, the author was supported by the Algebraic Combinatorics in Europe network, grant RTN2-2001-00059}
\maketitle
\begin{small}
\tableofcontents
\end{small}

\section{Introduction}

Let $\H$ be an affine Hecke algebra of type $B$ or $C$ with
possibly unequal labels (see \ref{aha} for its precise
definition).
%associated to the following data: a root system $R_0$, a lattice $X$ which is either the root lattice of $R_0$ or its weight lattice, and a length-multiplicative function $q: W \to \R_{>0}$ on the affine Weyl group $W=W_0 \ltimes X$ (where $W_0 \subset GL(X)$ is the finite Weyl group of $R_0$).
In this article we study the reducibility of the induction of a
discrete series representations of a parabolic subalgebra of $\H$.
We restrict ourselves to the case where the irreducible components of the induced representation have real central character. In particular, we may assume that the representation which is being induced has real
central character.  Under this assumption the results for type $C$
can be derived from those of type $B$ (see Remark \ref{BC}), so we
only compute the type $B$ case.

We study this reducibility by computing the relevant $R$-group. This is possible thanks to the recent extension of the
$R$-group theory by Delorme and Opdam to the case of affine Hecke
algebras. In the case of real reductive groups, this theory was
developed by Knapp, Stein and Harish-Chandra. The extension to
$p$-adic groups was done by Silberger. In each of the cases, one
constructs, starting from a subgroup of the Weyl group which is
computed using the Plancherel measure of the group (or the
algebra), intertwining operators for a certain parabolically induced
representation, and prove that these operators span its
centralizer algebra. The $R$-group consists of those elements
which yield nontrivial intertwiners. In general, the map $r
\mapsto E(r)$ from the $R$-group to the intertwiners yields a
projective representation.

We prove that in our situation, the $R$-group is always isomorphic to a certain $\Z_2^d$, where we can calculate $d$ explicitly using combinatorics involving the central character of the induced representation and the structure of the parabolic root subsystem which defines the Hecke subalgebra from which we induce.  This allows us to prove that the induced representation under consideration always decomposes into $2^d$ irreducible and pairwise inequivalent irreducible components. Moreover, we also prove that the a priori projective representation $r \mapsto E(r)$ is linear in our setting.

Suppose that $R_0$ is of type $B_n$ and let $q_1$, resp. $q_2$ be
the parameters of $\H$ corresponding to a simple reflection in a
long, resp. short root (see \ref{aha} and \eqref{q1q2}). Since
$q_1,q_2 \in \R_{> 0}$, there exists $m$ such that $q_2=q_1^m$.
The representation theory of $\H$ is largely determined by $m$ in
the sense that raising $q_1,q_2$ to a common power $q_1^\e,q_2^\e$
does not change the parametrization of the irreducible modules of
the corresponding Hecke algebras, even though the algebras
themselves are not isomorphic.

For all $m$ such that reducibility of induced discrete series representation in our setting turns out to occur, the Hecke algebra is isomorphic to the centralizer algebra of a certain parahorically induced representation of a reductive $p$-adic group $\G$ (see below). %In this setting, the corresponding problem of reducibility is very classical, since it arises naturally in the classification problem of the admissable dual of $\G$.
Thus, our results fit naturally into the theory of $R$-groups for
reductive $p$-adic groups of classical type. Goldberg has given an
explicit description of the $R$-groups of the groups of classical
Lie type $SO(n,F)$ and $Sp(2n,F)$ in \cite{goldberg}. He has
proven that for these groups, the $R$-group is isomorphic to
$\Z_2^d$ for some $d$. However, $d$ can not be calculated
explicity, other than by assuming knowledge of the irreducibility
of other parabolically induced discrete series representations
(those induced from a so-called basic parabolic subgroup). %Our approach of using the Hecke algebra rather than the $p$-adic group allows us to combinatorially express $d$ in terms of the central character of the induced discrete series representation, and the structure of the parabolic root system. Moreover, in this setting we can also show that the irreducible summands of the induced representation are in bijection with the irreducible characters of the $R$-group.

The outline of this paper is as follows.  In section 2 we review
some necessary material on affine and graded Hecke algebras of
arbitrary type. We review the relation between the affine Hecke
algebra and the representation theory of reductive $p$-adic
groups, explaining how the results in this paper connect to the
ones of Goldberg. In section 3 we review the portion of the theory
of the $R$-group developed by Delorme and Opdam which is relevant
for our purposes. Since our result fits naturally into a
combinatorial framework called generalized Springer correspondence
which we have introduced in \cite{S}, we quickly review some of
the involved constructions in section 4. The result of this paper
gives a partial affirmative answer to a conjecture in \cite{S}
which we also mention here. In section 5 we compute the $R$-group,
making use of these combinatorics, and prove the results mentioned above. For the Iwahori-Hecke case, we also give an interpretation of the $R$-group in terms of the component groups which come up in the Springer correspondence. This is done using the combinatorics recalled in section 4.

\section{The affine Hecke algebra}\

\subsection{Definition}\label{aha}Let ${\mathcal R}=(R_0,X,\check{R}_0,Y,\Pi_0)$ be a root datum.
By this we mean that $X$ and $Y$ are free abelian groups with a
perfect pairing $\langle \cdot,\cdot \rangle$ over $\Z$, that $R_0
\subset X$ is a reduced root system, that $\check{R}_0$ is the
dual root system of coroots of $R_0$, and that
$\Pi_0=\{\a_1,\dots,\a_n\}$ is a choice of simple roots of $R_0$.
We denote the corresponding set of positive roots by $R_0^+$.
Define $Q=\Z R_0 \subset X$, the root lattice. Let $W_0$ be the
Weyl group of $R_0$, that is, the finite reflection group
generated by the reflections $s_i: X \to X: x \to x-\langle x,
\ca_i \rangle \a_i$. Given $\Pi_L \subset \Pi_0$, we denote the
corresponding sub root system of $R_0$ by $R_L$, and
$W_L=W_0(R_L)\subset W_0$.

The extended affine Weyl group is defined as $W=W_0 \ltimes X$. It
admits a decomposition $W=W^a \rtimes \Omega$, with $\Omega \simeq
X/Q$ and $W^a=W_0 \ltimes Q$. This group $W^a$ is called the
affine Weyl group and is, like $W_0$, a Coxeter group. Let $S^m$
be the set of maximal coroots in $R_0^\vee$ (with respect to the
dominance ordering on $Y$). Then a set of Coxeter generators for
$W^a$ is given by $S=\{s_1,\dots,s_n\} \cup \{ s_\theta \mid
\theta^\vee \in S^m\}$ where $s_1,\dots,s_n$ are as above and
$s_\theta(x)=x-(-\langle x,\check{\theta} \rangle +1)(-\theta)$.
If $R_0$ is irreducible then $S^m$ consists of a single element
and we write $s_0=s_\theta$ for the corresponding reflection.
Being a Coxeter group, $W^a$ has a length function $l$. We can
uniquely extend it to a function on $W$, by defining $l(\omega)=0$
for $\omega \in \Omega$.

Choose a function $q: W \to \R_{>0}$ such that $q(ww')=q(w)q(w')$
whenever $l(ww')=l(w)+l(w')$. Such a function is determined by its
values on the simple reflections. Equivalently, it can be
characterized by a set of ``root labels" $q_{\ca} \in \R_{>0}$,
for $\a \in R_{nr}=R_0 \cup \{2\a \mid \ca \in 2Y \cap
\check{R}_0\}$, with the condition that $q_{\ca}$ depends only on
the $W_0$-orbit of $\a$. For $\a$ such that $2\a \notin R_{nr}$,
we formally define $q_{\ca/2}=1$. For an irreducible root system,
we define on simple affine roots $q(s_i)=q_{\ca_i}q_{\ca_i/2}$
($i=1,\dots,n$) and $q(s_0)=q_{\check{\theta}}$.

The affine Hecke algebra $\H=\H({\mathcal R},q)$ is then the
unique associative, unital complex algebra with $\C$-basis $T_w, w
\in W$ such that the multiplication satisfies
\[ \begin{array}{l} T_w T_{w'}=T_{ww'} \mbox{ if } l(ww')=l(w)+l(w');
 \\ (T_{s_i}+1)(T_{s_i}-q(s_i))=0 \mbox{ for }
 i=0,1,\dots,n.\end{array}\]

%\begin{rem}For certain $q:W\to \R_{>0}$, $\H(\mathcal{R},q)$ arises naturally in the representation theory of reductive $p$-adic groups, as the centralizer algebra of a certain induced representation (see \cite{lusuni1}, \cite{Morris}). In this case, one obtains a Plancherel measure preserving bijection between the irreducible representations of $\H$ and a certain subset of the irreducible admissible representations of the $p$-adic group. \end{rem}
\subsection{Irreducible tempered representations of $\H$}

\subsubsection{Bernstein-Zelevinsky presentation} Let $X^+=\{ x \in X \mid
\langle x,\ca \rangle \geq 0 \mbox{ for all } \a \in R_0^+\}$. For
$x \in X^+$, let $t_x$ denote the corresponding element of $W$ and
define $\theta_x=q(t_x)^{-1/2}T_{t_x}$. If $z \in X$, write it as
$z=x-y$ with $x,y \in X^+$, and put
$\theta_z=\theta_x\theta_y^{-1}$ (notice that all $T_w$ are
invertible). Let $\A$ be the subalgebra of $\H$ which has the
$\theta_x, x \in X$ as $\C$-basis. The action of $W_0$ on $X$ also
provides an action on ${\mathcal A}$. Let $\H_0$ be the subalgebra
of $\H$ with basis $T_w, w \in W_0$. According to what is known as
the Bernstein-Zelevinsky presentation, one has $\H=\H_0 \otimes \A
= \A \otimes \H_0$. Moreover Bernstein (unpublished, see
\cite{lusztig}), showed that the center $\cent$ of $\H$ is equal
to $\cent=\A^{W_0}$.

\subsubsection{Central character} It is well known that all irreducible representations
of $\H$ are of finite dimension. This follows from the fact that
$\H$ is finitely generated over $\cent$, and that by Dixmier's
version of Schur's Lemma, $\mathcal{Z}$ acts by scalars in an
irreducible representation $(\pi,V_\pi)$. Let
$T=\mbox{Hom}_\Z(X,\C^\times)=\mbox{Spec}(\A)$.  In an irreducible
representation $(\pi,V_\pi)$ of $\H$, $\cent$ acts by $W_0t \in
W_0\backslash T\simeq \mbox{Spec}(\cent)$. We call $W_0t$ (or, by
abuse of terminology, any of its elements) the central character
of $(\pi,V_\pi)$. Let $T$ have polar decomposition $T=T_uT_{rs}=
{\rm Hom}(X,S^1){\rm Hom}(X,\R_{>0})$. If $W_0t \subset T_{rs}$,
then we say that $(\pi,V_\pi)$ has real central character.

\subsubsection{Tempered; discrete series} Let $(\pi,V_\pi)$ be an
irreducible representation of $\H$. If, for $t \in T$, we have
$V_\pi(t)=\{ v \in V_\pi \mid (\pi(a)-t(a))v=0 \mbox{ for all } a
\in \A\} \neq 0$, then we say that $t$ is an $\A$-weight of
$(\pi,V_\pi)$. If every $\A$-weight $t$ of $(\pi,V_\pi)$ satisfies
$|t(x)| \leq 1$ for all $x \in X^+$, then we say that
$(\pi,V_\pi)$ is a tempered representation. If moreover $|t(x)|<1$
for all $0 \neq x \in X^+$, then we say that $(\pi,V_\pi)$ is a
discrete series representation. Let $\hat{\H}^t$ be the set of
equivalence classes of irreducible tempered representations of
$\H$. The theory of the Plancherel formula for the affine Hecke
algebra of \cite{O1} implies that the representations in
$\hat{\H}^t$ are precisely those which occur in the spectral
decomposition of the natural trace of $\H$ (cf. \cite{O1,DO}). We
will be concerned with $\Hrcc$, the representations in
$\hat{\H}^t$ with real central character. We denote by
$\hat{\H}^{ds}_\R \subset \hat{\H}^t_\R$ the subset of discrete
series representations.

\subsubsection{Residual cosets} We need some more terminology from
\cite{O1}. Let $L=r_LT^L$ be a coset of a subtorus $T^L$ of $T$.
Then we define a parabolic root system $R_L=\{\a \in R_0 \mid
\a(T^L)=1\}$. Let $R_L^p=\{\a \in R_L \mid
\a(r_L)=q_{\ca}q_{\ca/2}^{1/2} \mbox{ or }
\a(r_L)=-q_{\ca/2}^{1/2} \}$ and $R_L^z=\{\a \in R_L \mid
\a(r_L)=\pm 1\}$. Then we put $i_L=|R_L^p|-|R_L^z|$, and we say
that $L=r_LT^L$ is a {\it residual coset} if
$i_L=\mbox{codim}(T^L)$. Clearly the notion of residual coset is $W_0$-invariant. Usually we will choose $L$ in its orbit such that $R_L$ is a standard parabolic with simple roots $\Pi_L \subset R_L$.

If moreover $r_L$ can be chosen in
$T_{rs}$ then we say that $L$ is a {\it real} residual coset. A
residual coset $L=\{r_L\}$ (i.e., $T^L=\{1\}$) is called a residual point.  As
described in \cite[Ch. 7]{O1}, the classification of residual
points for graded Hecke algebras in \cite{HO} leads to a
classification of all residual cosets of $\H(\mathcal{R},q)$.

\subsubsection{$c$-function} For later reference, we also recall the following definition. Let $\a \in R_0$, then define $c_\a$ in the fraction field of $\mathcal{A}$ as
\begin{equation}
\label{cdef}
c_\a=\frac{(1+q_{\ca/2}^{-1/2}\theta_{-\a})(1-q_{\ca/2}^{-1/2}q_{\ca}^{-1}\theta_{-\a})}{(1+\theta_{-\a})(1-\theta_{-\a})}.
\end{equation}
We view $c_\a$ as a rational function on $T$. Notice that if $\ca \notin 2Y$, then by definition $q_{\ca/2}=1$, so $c_\a$ simplifies to $(1-q_{\ca}^{-1}\theta_{-\a})/(1-\theta_{-\a})$ in this case.

\subsubsection{Parabolic induction}\label{parabolic} Let $\Pi_L \subset \Pi_0$, then
${\mathcal R}^L=(R_L,X,\check{R}_L,Y,\Pi_L)$ is a root datum. Let
$X_L=X/(X \cap (R_L^\vee)^\perp)$ and $Y_L=Y \cap \Q R_L^\vee$.
Then ${\mathcal R}_L=(R_L,X_L,\check{R}_L,Y_L,\Pi_L)$ is also a
root datum. Restriction of $\{q_{\ca} \mid \a \in R_{nr}\}$ to
$R_{L,nr}\subset R_{nr}$ yields $q^L$ on $W_L \ltimes X$ and $q_L$
on $W_L \ltimes X_L$. Let $\H^L=\H({\mathcal R}^L,q^L)$ and
$\H_L=\H({\mathcal R}_L,q_L)$ be the associated affine Hecke
algebras. Let $T_L={\rm Hom}(X_L,\C^*)$. Define also the lattice
$X^L=X/(\R R_L \cap X)$ and the torus $T^L={\rm Hom}(X^L,\C^*)$.
It follows from the results in \cite{O1} that the representations
in $\hat{\H}^t$ all arise as irreducible components of
\begin{equation}
\label{indL} \pi(\Pi_L,\delta,t^L)=\mbox{Ind}_{\H^L}^\H(\delta \circ
\phi_{t^L}),
\end{equation}
where $\delta \in \hat{\H}_{L,\R}^{ds}$, $t^L \in T^L_u$, and
$\phi_{t^L}$ denotes the map $\H^L \to \H_L$ which takes $T_w
\mapsto T_w$ for $w \in W_L$, and sends $\theta_x \mapsto
t^L(x)\theta_{\bar{x}}$ with $x \mapsto \bar{x}$ the projection $X
\to X_L$. Moreover, a representation in $\Hrcc$ arises as an
irreducible summand in \eqref{indL} for a unique $(\Pi_L,\delta,1)$ if
we consider the tripels $(\Pi_L,\delta,1)$ up to a suitable
equivalence relation (which in the case of real central character
and a root system $R_0$ of type $B_n$ is just Weyl group conjugacy
of $L$, cf. Lemma. \ref{identity} below).

\subsubsection{Classification of central characters of $\Hrcc$}\label{cc}
In \cite{O1}, Opdam  obtains a classification of the central
characters of the irreducible tempered representations of $\H$,
hence in particular of the representations in $\Hrcc$. The central
characters of $\Hrcc$ are precisely the $W_0$-orbits of the
centers of the real residual cosets. This bijection is such that
the central character of $(\pi,V_\pi) \in \Hrcc$ which occurs in
the induced representation $\pi(\Pi_L,\delta,1_L)$ is the center of a
residual coset $M$ with $R_L=R_M$. Thus, the central character of
a discrete series representation is a residual point, and
conversely.

\subsubsection{The $R$-group}  Let $L=r_LT^L$ be a real residual coset.
In general, the description of $\Delta_{W_0r_L} \subset \Hrcc$ of
irreducible representations with central character $W_0r_L$ is
unknown, apart from the fact that it is a non-empty and finite
set. For a Hecke algebra $\H$ of type $A$, it is known that
$\hat{\H}^{ds}_\R$ consists of a single one-dimensional
representation. If $q_{\ca}>1$, the corresponding representation
of the symmetric group is the sign representation. The discrete
series representation itself is referred to as the Steinberg
character. For arbitrary root systems, to classify
$\Delta_{W_0r_L}$ for all central characters of $\Hrcc$, one needs
to classify all $\hat{\H}^{ds}_{L,\R}$, and understand the
induction \eqref{indL}. This latter problem is adressed in recent
work of Delorme and Opdam, \cite{DO2}, where they give a description of a
group which controls the decomposition into irreducible components
of \eqref{indL}. They call this group the $R$-group, in analogy
with the corresponding group for real or $p$-adic groups (see e.g.
\cite{arthur}).

\subsection{Hecke algebras and reductive $p$-adic
groups}

Let $\bf{G}$ be a connected split reductive group defined over a
non-archemidean non-discrete local field $F$ of characteristic
zero, with ring of integers $\O$ and residue field $\F_q$. Let $G$
be the group of $F$-rational points of $\bf{G}$.

%There are two important constructions which yield representations of $G$ that play a role in the parametrization of its admissible dual.

\subsubsection{Parabolic induction in $G$} Let $P=MN$ be a parabolic subgroup with Levi
subgroup $M$ and unipotent radical $N$. Let $\rho$ be an
irreducible, admissible representation of $M$. Then we can extend
$\rho$ across $N$ to become a representation of $P$, and then
induce to obtain $I(P,\rho)={\rm Ind}_P^G(\rho)$. The compactness
of $G/P$ implies the admissibility of $I(P,\rho)$. It is known
that for each irreducible admissible representation $\pi$ of $G$,
there exists a parabolic subgroup $P=MN$ and an irreducible
supercuspidal representation $\rho$ of $M$ such that $\pi$ is a
subrepresentation of $I(P,\rho)$. Similarly, the $I(P,\rho)$ where
$\rho$ is a discrete series representation, yield the tempered
spectrum of $G$.

\subsubsection{Unipotent representations; the Hecke algebra} %Let $K=G(\O)$ be the maximal compact subgroup of $G$ which we obtain by taking the $\O$-rational points of\marginpar{ik stel voor deze definities bekend te veronderstellen en weg te laten} $\bf{G}$. Then the map $\O \to \F_q$ induces a map $K \to G(\F_q)$. Let $\I$ be the inverse image in $K$ of a Borel subgroup of $G(\F_q)$. Any subgroup of $G$ conjugate to $\I$ is called an Iwahori subgroup of $G$. Any open and compact subgroup of $G$ which contains an Iwahori subgroup is called a parahoric subgroup.
Let $\P \subset G$ be a parahoric subgroup, with pro-unipotent radical $U$.
Then the quotient $M=\P/U$ is isomorphic to a Levi subgroup of
$G(\F_q)$, hence in particular it is a finite reductive group of
Lie type. Let $\sigma$ be a cuspidal unipotent representation of
$M$, and let $\sigma^\vee$ be its contragredient. We inflate
$\sigma^\vee$ across $U$ to become a representation of $\P$, and
then form the (compactly) induced representation ${\rm
c-Ind}_\P^G(\sigma^\vee)$. Let $\H(G,\P,\sigma)={\rm End}_G({\rm
c-Ind}_\P^G(\sigma^\vee))$ be its centralizer algebra (the
notation is such as to obtain the equivalence of categories
below). By theorems of Morris, \cite{Morris} and Lusztig,
\cite{lusuni1}, $\H(G,\P,\sigma)$ is isomorphic to an (extended)
affine Hecke algebra associated to a root system whose rank is the
difference between the split ranks of $G$ and $M$, with explicitly
determined possibly unequal labels.

%We recall the following from \cite{Morris1999}. Fix a parahoric subgroup $\P$ of $G$, with pro-unipotent radical $U$ and let $(\sigma,V)$ be the inflation ${\rm inf}_M^\P(\sigma^\vee)$ of the contragredient of a cuspidal irreducible representation $(\sigma,V)$ of the reductive quotient $M=\P/U$ to $\P$. Remark that $\sigma^\vee$ viewed as a representation of $\P$ is still irreducible.

%%%%%Let $\H(G,\P,\sigma)={\rm End}_G(c-{\rm Ind}_\P^G(\sigma^\vee))$ be the centralizer algebra of the compactly induced representation of $(\sigma^\vee,V^\vee)$.

Let $\mathfrak{SR}_\sigma(G)$ be the category of smooth
representations of $G$ which are generated by their
$\sigma$-isotypic component. Then by \cite[Theorem 4.3]{BK1998},
there is an equivalence of categories between
$\mathfrak{SR}_\sigma(G)$ and the category $\H(G,\sigma)-{\rm
Mod}$ of modules of $\H(G,\P,\sigma)$.

%If $\sigma$ is a cuspidal unipotent representation of $M=\P/U$, then it is known (see \cite{Morris}) that $\H(G,\P,\sigma)$ is isomorphic to an affine Hecke algebra with possibly unequal labels, associated to a root system whose rank is the difference between the split ranks of $G$ and $M$.

The prototype is when one takes $\P=\I$, an Iwahori subgroup of
$G$ and $\sigma=1$. Suppose ${\bf G}$ has root datum
$(R_0^\vee,Y,R_0,X,\Pi_0^\vee)$. Then $\H(G,\I,1)$ is the
Iwahori-Hecke algebra with root datum $(R_0,X,R_0^\vee,Y,\Pi_0)$
and labels $q(s)=q$, the cardinality of
the residue field of $F$. %In particular, the torus $T={\rm Hom}(X,\C^*)$ is isomorphic to the torus of unramified characters of a split maximal torus ${\bf A}$ in $G$. Indeed, $A={\rm Hom}_\Z(Y,F^*), A(\O)={\rm Hom}_\Z(Y,\O^*)$ and so $A/A(\O) \simeq {\rm Hom}_\Z(Y,F^*/O^*)= {\rm Hom}_\Z(Y,\Z)=X$, and $T={\rm Hom}(X,\C^*)$.
The equivalence of categories ${\mathfrak S R}_\sigma(G) \to
\H(G,\I,\sigma)-{\rm Mod}$ is given by taking Iwahori-invariant
vectors, that is by the map $V \mapsto V^\I$, which is naturally a
$\H=\H(G,\I,1)$-module.
\subsubsection{Relation between parabolic inductions in $\H$ and $G$}\label{para}

For the general case of an affine Hecke algebra $\H(G,\P,\sigma)$,
it does not seem to be written up in the literature what the
precise relation is between parabolic induction in
$\H(G,\P,\sigma)$ and in $G$.

However in the Iwahori-spherical case where $\P=\I$ and $\sigma
=1$, the following is known from the work of Jantzen, cf.
\cite{jantzen95}.  Let $P=MN$ be a parabolic subgroup with
standard Levi subgroup $M$. Let $R_M$ be the root system of $M$.
Let $\rho$ be a representation of $M$ which is generated by its
$\I_M:=\I \cap M$-fixed vectors (note that $\I_M$ is an Iwahori
subgroup of $M$). Let $w_0$ (resp. $w_M$) be the longest element
of $W_0$ (resp. of $W_0(R_M)$) and let $w^M=w_0w_M$. Let
$M'=w_0Mw_0=w^MM(w^M)^{-1}$. Let $\H^{M'}$ be the subalgebra of
$\H$ associated to $R_{M'}$ as in \ref{parabolic}. Then, according
to \cite[Prop. 2.1.2]{jantzen95} (the case where $M$ is a maximal
split torus goes back to \cite{reeder}), we have
\begin{equation}
\label{jantzen} ({\rm Ind}_P^G(\rho))^\I \simeq \H
\otimes_{\H^{M'}} (w^M \cdot \rho^{\I_M}).
\end{equation}

This means that the results which we will describe, relevant to
the right hand side, can be interpreted as results on the
reducibility of the induced representation on the left hand side,
and thus refine the work of Goldberg in these cases.

%In the special case $\P=\I, \sigma=1$ then $\H(G,\I,1)$ is isomorphic to the (equal label) Iwahori-Hecke algebra, as shown by Iwahori and Matsumoto in \cite{IM}.

\subsubsection{Goldbergs results} For the induced representation $I(P,\rho)$, the decomposition into irreducibles of ${\rm Ind}_P^G(\rho)$ is determined by a finite group, referred to as the $R$-group. Roughly speaking, the dual of this group parametrizes the irreducible constituents. For $G=GL(n,F)$, it is known that all $R$-groups are trivial. For the split groups of classical type $SO(n,F),Sp(2n,F)$, the $R$-group has been computed by Goldberg in \cite{goldberg}. He shows that if $\rho$ is a discrete series representation, then its $R$-group $R(\rho) \simeq \Z_2^d$ for some $d$. We discuss this in some more detail for the case where $G=SO(2n+1,F)$, the case where $G=Sp(2n,F)$ being analogous (we are not concerned with the type $D$ case where $G=SO(2n,F)$ in this paper). Let $M$ be a Levi subgroup of $G$, then $M$ can be written as
\begin{equation}
\label{levi} M \simeq GL(m_1)^{n_1} \times \dots \times
GL(m_r)^{n_r} \times SO(2m+1),
\end{equation}
where $\sum_im_in_i+m=n$ and the $m_i$ are all different. If $r=1$
then $M$ is called a basic Levi subgroup, and $P=MN$ a basic
parabolic. Given $M$ as in \eqref{levi} and a discrete series
representation $\rho$ of $M$, then $\rho=\rho_1 \otimes \dots
\otimes \rho_t \otimes \rho_B$, where $\rho_i$ is a discrete
series representation of $GL(m_i)^{n_i}$ and $\rho_B$ is a
discrete series representation of $SO(2m+1)$. Let $R_i$ be the
$R$-group of the representation $\rho_i \otimes \rho_B$, w.r.t to
the induction to $SO(2(m_in_i+m)+1)$. Then Goldberg has shown that
\begin{equation}
\label{Rprod} R(\rho) \simeq R_1 \times \dots \times R_t.
\end{equation}
 Every $R_i$ is of the form $R_i \simeq \Z_2^{d_i}$. Let $\rho_i=\rho_{i,1} \otimes \dots \otimes \rho_{i,n_i}$ where each $\rho_{i,j}$ is a discrete series representation of $GL(m_i)$. Then, if $P_i'=M_i'N_i'$ is the parabolic subgroup of $G_i \simeq SO(2(m_i+m)+1)$ such that $M_i'\simeq GL(m_i) \times SO(2m+1)$, we have
\begin{equation}
\label{di} d_i=\sharp\{ [\rho_{i,j}] \mid {\rm
Ind}_{P'_i}^{G_i'}(\rho_{i,j}\otimes \rho_B) {\rm \ is \
reducible}\}.
\end{equation}
Here, $[\rho_{i,j}]$ denotes the equivalence class of
$\rho_{i,j}$. In \cite{MT2002}, M\oe glin and Tadi\'c compute the
numbers $d_{i,j}$ in terms of so-called Jordan pairs. However,
these pairs seem not to be easy to compute.

%\subsubsection{}
For affine Hecke algebras, Delorme and Opdam have shown in
\cite{DO2} that if one considers unitary parabolic induction of
discrete series representations of a parabolic quotient algebra,
one can also construct the $R$-group, and that it has the same
significance. In this article, we explicitly compute this
$R$-group in the case where the induced representation has real
central character. We will show that in this special case the
analogues of \eqref{Rprod} and \eqref{di} hold. However for the
analogue of \eqref{di}, the assumption of real central character
implies that each $\rho_{i,j}$ is the Steinberg representation, so
we are in the simpler situation where $d_i \in \{0,1\}$. We
compute $d_i$ combinatorially, using a certain Young tableau which
is defined in terms of the parameters of the Hecke algebra, the
type of the parabolic and the central
character of the induced discrete series representation. %In particular the $R$-group, hence the reducibility of $I(P,\rho)$, depends only on the central character of $\rho$, rather than on $\rho$ itself.

\subsection{The graded Hecke algebra}\label{GHA}
In the proof of Theorem \ref{deco} which describes the
decomposition into irreducibles of the representations
$\pi(\Pi_L,\delta,1_L)$ for $\delta \in \hat{\H}^{ds}_{L,\R}$, we will
need the graded Hecke algebra and its relation to the affine Hecke
algebra, as described by Lusztig in \cite{lusztig}.

Let $\a \in R_0$. Then we put
\begin{equation}
\label{kq}
k_\a= \mbox{log}(q_{\ca}q_{\ca/2}^{1/2}).
\end{equation}
Notice that if $\ca \notin 2Y$, we simply have $k_\a={\rm log}(q_{\ca})$.
Let $\real^*=X \otimes_\Z\R$, and define ${\mathcal R}^{deg}=(R_0,\real^*,R_0^\vee,\real,\Pi_0)$ (we call this the degenerate root datum associated to ${\mathcal R}$). Then the graded Hecke algebra associated to ${\mathcal R}^{deg}$ and the label function $k: \a \mapsto k_\a$ is defined to be
\[ \Hgr=\Hgr({\mathcal R}^{deg},k)=\C[W_0] \otimes S(\real^*_\C),\]
where we have the relations
\[ x\cdot s_\a - s_\a \cdot s_\a(x)=k_\a \langle x, \ca \rangle,\]
for $x \in \real^*_\C, \a$ simple.
As reviewed in \cite{S}, it follows from theorems of Lusztig in \cite{lusztig} that there exists a natural bijection between $\Hrcc$ and $\Hgrrcc$ (defined analogously). This bijection is such that if $(\pi,V) \in \Hgrrcc$ has central character $W_0\ga$, then the corresponding module in $\Hrcc$ also has $V$ as underlying vector space, and its central character is $W_0{\rm exp}(\ga)$ where ${\rm exp}: Y\otimes_\Z\R \to T_{rs}$ is the exponential map (we have $\ga \in Y\otimes_\Z\R$ since we assume it is a real central character).

For the graded Hecke algebra, one has the analogous notion of residual subspace. The real residual subspaces are the logarithms of the real residual cosets in $T$. As in the affine case, the central characters of $\Hgrrcc$ are the $W_0$-orbits of the centers of the residual subspaces.

\begin{rem}\label{BC} Suppose that $\H=\H({\mathcal R'},q')$ where
${\mathcal R}=(R_0',X',R_0',Y',\Pi_0')$ and $R_0'$ has type $C_n$.
Then, by the above mentioned theorems of Lusztig, there is a
natural bijection $\H'^t_{\R} \leftrightarrow \Hgr'^t_\R$, where
$\Hgr'$ is the graded Hecke algebra associated to ${\mathcal
R}'^{deg}$ and $k'$ depending on $q'$ as in \eqref{kq}. Let $k_1'$
(resp. $k_2'$) be the labels of the roots of type $\pm e_i \pm
e_j$ (resp. $\pm 2e_i$), expressed in standard basis vectors of
$X' \otimes_\Z \R$. Let $R_0$ be the root system of type $B$
obtained from $R_0'$ by replacing the $\pm 2 e_i$ with $\pm e_i$.
Let $k_1=k_{\pm e_i \pm e_j}=k_1'$, and $k_2=k_{\pm e_i}=
\frac{1}{2}k_2'$. It is not hard to see that $\Hgr' \simeq \Hgr$,
where $\Hgr$ is attached to the data ${\mathcal R}^{deg}$,
obtained from ${\mathcal R}^{deg}$ by replacing $R_0'$ with $R_0$,
and label function $k$. Thus, the problem of reducibility of an
induced representation with real central character for an affine
Hecke algebra of type $C$ is equivalent to a corresponding problem
for an affine Hecke algebra of type $B$. Therefore we will only
treat the type-$B$ case.
\end{rem}

\section{The $R$-group}
We now review the definition and properties of the $R$-group,
following \cite{DO2}.  Let $\Xi$ be the (cf. \cite{O1}) set of
objects of the  ``groupoid of standard induction data'', that is,
an element of $\Xi$ is a triple $\xi=(\Pi_L,\delta,t^L)$, where
$\Pi_L \subset \Pi_0$ defines a standard parabolic, $\delta \in
\hat{\H}^{ds}_L$, and $t^L \in T^L$. If moreover $t^L \in T^L_u$,
then we write $\xi \in \Xi_u$.

Suppose that $\delta$ has central character $W_Lr_L \subset T_L$,
then we will say that $\xi$ has central character $W_0r_Lt^L$. Let
\[ \pi(\xi)=\pi(\Pi_L,\delta,t^L)=\mbox{Ind}_{\H_L}^\H(\delta \circ \phi_{t^L}). \]
The central character of the irreducible components of $\pi(\xi)$
is $W_0r_Lt^L \subset T$. Define
\[ {\mathcal W}_{L,L}=\{k \times w \in  (T^L\cap T_L) \times W_0 \mid w(L)=L\},\]
and (cf. \cite[p. 34]{DO2})
\begin{equation}
\label{Wxixi} {\mathcal W}_{\xi,\xi}=\{ g= k \times w \in
{\mathcal W}_{L,L} \mid g\cdot \xi=(\Pi_L,\Psi_g(\delta),kw(t^L))=\xi\}.
\end{equation}

This means the following, cf. \cite{O1}. If $w(L)=L$, then $w$
induces an automorphism of affine Hecke algebras $\psi_w : \H_L
\to \H_L$. If $k \in K_L=T^L \cap T_L$ then we also have an
automorphism $\psi_k: \H_L \to\H_L: \theta_{\bar{x}}N_w\mapsto
k(\bar{x})\theta_{\bar{x}}N_w$ for $\bar{x} \in X_L$ and $w \in
W_L$. Let $g=k \times w$, and define $\psi_g: \H_L \to \H_L$ to be
the composition $\psi_g=\psi_k \circ \psi_w$. This leads to a
bijection $\Psi_g: \Delta_{W_Lr_L} \to \Delta_{k^{-1}W_Lw(r_L)}$
given by $\Psi_g(\delta) \simeq \delta \circ \psi_g^{-1}$.

\begin{lemma} \label{identity}If $R_0$ is of type $B_n$ and $\xi=(\Pi_L,\delta,1)$ has real central character, then $\Psi_g$ is the identity map for all $g \in {\mathcal W}_{\xi,\xi}$.
\end{lemma}

\pf:  The fact that $t^L=1$ implies $k=1$ for all $g=k \times w
\in {\mathcal W}_{\xi,\xi}$. Suppose that the root system $R_L$ is
of type $\prod_{i} A_{\l_i} \times B_l$ for some $\l_i,l$. Then
$w(L)=L$ implies that $w$ acts as the identity on the part of type
$B_l$. On the parts of type $A$, it acts as a combination of
diagram automorphisms and permutations of $A$-factors of equal
rank. Since $W_Lr_L$ is the central character of a discrete series
representation of $\H_L$, it is the product of the central
characters of discrete series representations of affine Hecke
algebras of type $A_{\l_i}$ and a central character of discrete
series representations of an affine Hecke algebra of type $B_l$.
An affine Hecke algebra of type $A$ which occurs as a factor in
$\H_L$ has only one discrete series representation with real
central character (the Steinberg representation). Therefore
$\Psi_g$ acts trivially on $\Delta_{W_Lr_L}$.\qed

It follows that, if $R_0$ has type $B_n$ and for $\xi$ with real
central character, we have
\begin{equation}
\label{rcc}
{\mathcal W}_{\xi,\xi}= \{ (w,1) \mid w(L)=L\}.
\end{equation}
In particular, ${\mathcal W}_{\xi,\xi}$ depends only on $L$ and not on $\delta$.

Let $\t^*=X\otimes_\Z\C \supset R_0$. For $L \subset I$, define
\[ \t^{L,*}=\{ \l \in \t^* \mid \l(\ca)=0 \mbox{ for all } \a \in
L\}.\]

\begin{defn}(cf. \cite[Definition 5.1]{DO2}) Let $\a \in R_0 \backslash R_L$, then we denote by $[\a]^L \in \t^{L,*}$ the restriction of $\a \in \t^*$ to $\t^{L,*}$. Let $\beta = [\a]^L$
for $\a \in R_0\backslash R_L$. Let $\xi=(\Pi_L,\delta,t^L) \in \Xi$
where $\delta$ has central character $W_Lr_L$. Then we define
\begin{equation}
\label{cbeta}
 c_\beta(\xi)=\prod_{\{\a\mid [\a]^L \in \R_{>0}\beta\}}c_\a(r_Lt^L).
\end{equation}
%and \[ c(\xi)=\prod_{\beta}c_\beta(\xi).\]
\end{defn}

\begin{lemma}
With the above notation, $c_\beta(\xi)$ does not depend on the
choice of $r_L$ in its $W_L$-orbit.
\end{lemma}

\pf: %Recall $c_\beta(\xi)=\prod_{\{\a'\mid [\a']^L \in \R_{>0}\beta\}} c_{\a'}(r_L)$.
Let $w \in W_L$. Then, for any root $\a$ we have $\a= \a_{\mid
\t_L} + \a_{\mid \t^L}$. Since $R_L \subset \t_L$, we have
$w(\a)=\a_{\mid\t^L}+w(\a)_{\mid \t_L}$, hence $[w(\a)]^L=[\a]^L$.
The lemma follows.\qed

\begin{defn} Let $\beta=[\a]^L$ for $\a \in R_0-R_L$. Then we say that $\beta$ is primitive, if for every $\a' \in R_0-R_L$ we have $[\a']^L \in \R_{>0}\beta \Rightarrow [\a'] \in \Z_{>0}\beta$.
\end{defn}

\begin{thm} (\cite[Prop. 6.5]{DO2}) Let $\xi\in \Xi_u$. Let $\beta=[\a]^L$ be primitive. Then the pole order of \eqref{cbeta} is at most
one. Moreover, the primitive $\beta=[\a]^L$ for which \eqref{cbeta}
has a pole form an integral root system in $\t^{L^*}$. We will
denote this root system by $R_0(\xi)$.
\end{thm}

A Weyl group element $w \in W_0$ for which $w(L)=L$ acts on the restricted roots by $w([\a]^L)=[w(\a)]^L$. Therefore one can define (cf. \cite[Def. 6.8]{DO2})
\[ R(\xi)=\{g \in {\mathcal W}_{\xi,\xi} \mid
g(R_{0,+}(\xi))=R_{0,+}(\xi)\}.\] Suppose that $\xi \in \Xi_u$.
Then, according to \cite[Prop. 6.7]{DO2}:
\[ {\mathcal W}_{\xi,\xi} \simeq W_0(\xi) \rtimes R(\xi),\]
where $W_0(\xi)$ denotes the Weyl group of $R_0(\xi)$.

For every $w \in {\mathcal W}_{\xi,\xi}$, Opdam has defined (cf. \cite[4.4]{O1}) an intertwining operator $E(w,\xi)$ of $\pi(\xi)$.
\begin{thm} \cite[Theorem 6.3]{DO2} Let $\xi \in \Xi_u$.
If $w \in W_0(\xi)$ then $E(w,\xi)=\l Id_{\pi(\xi)}$.
\end{thm}
Moreover, the elements of $R(\xi)$ lead to the nontrivial $\H$-endomorphisms of
$\pi(\xi)$ and determine its decomposition into irreducibles:

\begin{thm} \label{DOthm}(\cite[Theorem 7.4]{DO2}) For all $\xi\in\Xi_u$, one has
\begin{itemize}
\item ${\rm End}_\H(\pi(\xi))=\sum_{r \in R(\xi)}\C E(\xi,r)$;
\item ${\rm dim}({\rm End_\H}(\pi(\xi))=|R(\xi)|$;
\item $E(r,\xi)E(r',\xi)=\eta_\xi(r,r')E(rr',\xi)$ for all $r,r' \in R(\xi)$, where $\eta_\xi$ is a 2-cocycle.% $\eta_\xi(r_1,r_2)\eta_\xi(r_1r_2,r_3)=\eta_\xi(r_1,r_2r_3)\eta_\xi(r_2,r_3)$ for all $r_1,r_2,r_3$ in $R(\xi)$.
\end{itemize}
\end{thm}

%\subsubsection{Reducibility of $\pi(\xi)$}\label{decomp}
The first property says that the intertwining operators $E(r,\xi)$
with $r \in R(\xi)$ span the space of intertwiners of $\pi(\xi)$,
and by the second property they are linearly independent, hence
form a basis of this space.

If $\eta_\xi$ splits, i.e., satisfies
\begin{equation}
\label{split}
 \eta_\xi(r,r')=\sigma_\xi(rr')\sigma_\xi(r)^{-1}\sigma_\xi(r')^{-1}
\end{equation}
for some $\sigma_\xi: R(\xi) \to \C^*$, then by standard arguments (see \cite{arthur}) one obtains from Theorem \ref{DOthm} a bijection between the irreducible representations of $R(\xi)$ and the irreducible components of $\pi(\xi)$. %
%Equation \eqref{split} is equivalent to $\bar{\eta}_\xi \in H^2(R(\xi),\C^*)$ being trivial.
We will see below that if $R_0$ is of type $B_n$ and $\xi \in
\Xi_u$ has real central character, then we are in this situation.

\section{Conjecture on $\Hrcc$ for $R_0$ of type $B$}\label{combi}

Given a orthonormal basis $e_1,\dots,e_n$ of $Q \otimes_\Z \R$, we
let $R_0=\{\pm e_i \pm e_j \mid 1\leq i,j \leq n\} \cup \{\pm
e_i\mid 1\leq i \leq n\}$ and we choose as simple roots
$\a_1,\dots,\a_n$ where $\a_i=e_i-e_{i+1}$ for $1 \leq i <n $ and
$\a_n=e_n$. One has $\theta=e_1$. Define (cf. \eqref{kq}, noticing
that $\ca_1 \notin 2Y$ for any choice of root datum),
\begin{equation}
\label{q1q2}
 q_1=q_{\ca_1}, q_2=q_{\ca_n}q_{\ca_n/2}^{1/2}.
\end{equation}
We will always assume that $q_1 \neq 1$. Furthermore we assume that
$\mbox{rank}(X)=\mbox{rank}(Q)$, a necessary condition in order to
have discrete series representations. Thus, $X$ is either the root lattice $Q$ or the weight lattice. In the latter case $Y=\check{Q}$ is the coroot lattice, hence $R_{nr}=R_0$ and all $q_{\ca/2}=1$.

\subsection{Residual points} First we recall some facts from \cite{HO} and
\cite{S} (see also \cite{proefschrift}).  Define (since we assume
that $q_1 \neq 1$) $m$ by $q_2=q_1^m$. We say that the parameters
$q_1,q_2$ are {\it generic} if $m \notin
\pm\{0,1/2,1,3/2,\dots,n-1\}$. Otherwise, we call them {\it
special}. According to \cite{HO}, the real residual points for
$\H$ are, for generic parameters, indexed by the partitions of
$n$. Given $\l \vdash n$, let $Y(\l)$ be its Young tableau. For a
square $\Box \in Y(\l)$, let $c(\Box)$ be its content (its column
coordinate minus its row coordinate). Then we define
$c(\l;q_1,q_2)$ to be the point with coordinates
$q_1^{c(\Box)}q_2$. We do not specify an order on the coordinates,
since we are interested in $W_0$-orbits only.

\begin{figure}[htb]
\begin{center}
{\includegraphics[angle=0,scale=0.38]{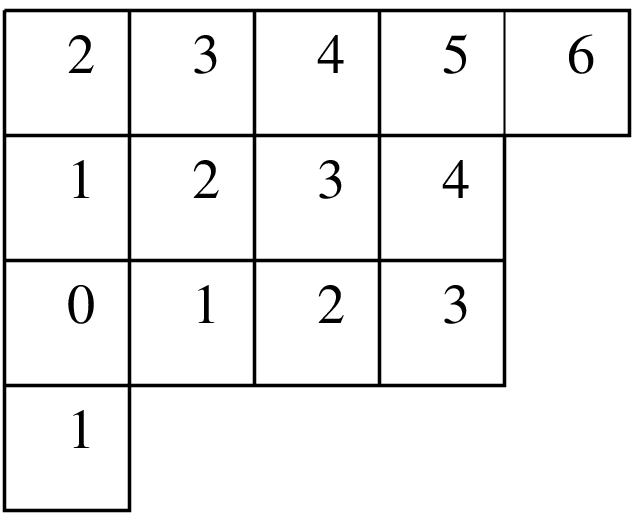}}
\end{center}
\caption{$T_2(\l)$ for $\l=(1445)$.}
\label{mtab}
\end{figure}

 For $m \in \R$, we define for a partition $\l$ its $m$-tableau $T_m(\l)$ to
be its Young tableau where box $\Box$ is filled with entry
$e_m(\Box)=|c(\Box)+m|$. An example for $m=2$ is given in Figure \ref{mtab}.
Then we recall from \cite{S} the
splitting map $\S_m$:

\begin{defn} Let $\l \vdash n$. The splitting
procedure $\S_m$ divides $T_m(\l)$ into horizontal and vertical
blocks, by subsequently enclosing the biggest remaining entry in
the thus far unselected part into a horizontal or vertical box, in
such a way that at every moment the thus far created boxes form
the diagram of a partition.
\end{defn}

\begin{figure}[h]
{\includegraphics[width=10cm]{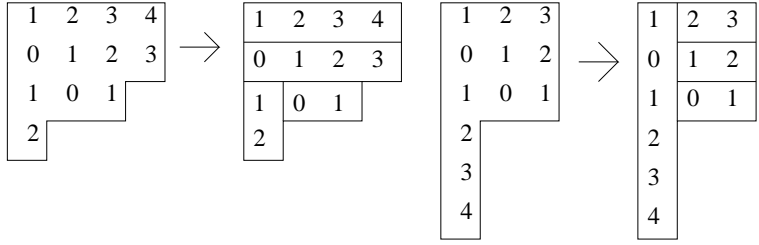}}
\caption{Examples of the splitting map}
\label{fig:demosplits}
\end{figure}

Figure \ref{fig:demosplits} contains, for $m=1$, two examples of $\l$ such that $\S_m$ is well-defined on $Y(\l)$.

The map $\S_m$ is not well-defined for all $m$ and all $\l$: it
may happen that the maximal remaining
entry occurs twice. But we have:

\begin{prop}\cite{S}
Let $q_1>0,q_2>0$ and $q_2=q_1^m$. Then $c(\l;q_1,q_2)$ is a
residual point if and only if $\S_m$ is well-defined on $T_m(\l)$.
\end{prop}

Let $\P(n,2)$ be the set of bipartitions $(\xi,\eta)$ of total
weight $n$. If $\S_m$ is well-defined on $\l$, then we define
$\S_m(\l)=(\xi,\eta) \in \P(n,2)$, where $\xi$ consists of the
lengths of the horizontal blocks, and $\eta$ of the lengths of the
vertical blocks in the split tableau. A $ 1 \times 1$-block $\Box$
is considered to be horizontal, resp. vertical, if $c(\Box)>-m$,
resp. $c(\Box)<-m$ (that is, if $\Box$ lies above, resp. below,
the zero-diagonal).

If $m \notin \pm\{0,1/2,1,\dots,n-1\}$ then $\S_m$ is well-defined
for all $\l$. For special parameters, one may have
$W_0c(\l;q_1,q_2)=W_0c(\mu;q_1,q_2)$ even if $\l \neq \mu$ and
this point need not be residual.

\subsection{Symbols and Springer correspondents} Conjecture \ref{con} below is formulated using
a generalization of the symbols which were introduced by Lusztig
in \cite{luscells} to describe the Springer correspondence. These
symbols are defined as follows. Here and in the sequel, it is our convention to write the
parts of a partition in increasing order: $\l=(0\leq \l_1 \leq \l_2 \leq \dots )=(0^{r_0}1^{r_1}2^{r_2}\dots)$.

Let $m \in \Z$. Then the $m$-symbol $\Lambda^m(\xi,\eta)$ of $(\xi,\eta) \in \P(n,2)$ is the two-lined array whose top line
contains $\xi_1,\xi_2+2,\xi_3+4,\dots$ and whose bottom line
contains $\eta_1,\eta_2+2,\eta_3+4,\dots$, where we assume that
$l(\xi)=l(\eta)+m$ by adding zeroes to $\xi$ or $\eta$, but such that $\eta_1$ or $\xi_1$ is non-zero. In
the symbol we write the entry of $\xi_i$ to the left of the
corresponding entry of $\eta_i$ if $m>0$, and to the right of the
corresponding entry of $\eta_i$ if $m<0$. For $m=0$, we have two
symbols, the $+0$-symbol and the $-0$-symbol. For example,
$\Lambda^2(12,3)$, $\Lambda^{+0}(12,3)$, $\Lambda^{-0}(12,3)$ and
$\Lambda^{-2}(12,3)$ are, respectively
\[  \left( \begin{array}{lllll}0&&3&&6 \\ &1&&& \end{array}
\right),\left(\begin{array}{llll}1&&4& \\ &0&&5 \end{array}
\right),\left(\begin{array}{lllll}&1&&4 \\ 0&&5& \end{array}
\right),\left(\begin{array}{lllllll}&1&&4&&&\\ 0&&2&&4&&9
\end{array} \right).\]
If $m \notin \Z, m \in \half\Z$, then the $m$-symbol
$\Lambda^m(\xi,\eta)$ is defined similarly. We ensure that
$l(\xi)=l(\eta)+sgn(m)(|m|+1/2)$. The symbol $\Lambda^m(\xi,\eta)$
has entries $\xi_1,\xi_2+2,\xi_3+4,\dots$ in the top row and
$\eta_1+1,\eta_2+3,\eta_3+5,\dots$ in the bottom row.

Suppose that $\Lambda^m(\xi,\eta)$ and $\Lambda^m(\a,\beta)$ have
the same entries with the same multiplicities. Then we say that
these symbols are {\it similar}, which we denote by $(\xi,\eta)
\sim_m (\a,\beta)$. The similarity class in $\P(n,2)$ under $\sim_m$ of
$(\xi,\eta)$ is denoted by $[(\xi,\eta)]_m$.

\subsubsection{Truncated induction}Let $(\xi,\eta)\in \P(n,2)$. Then the $m$-symbol
$\Lambda^m(\xi,\eta)$ can be written as
$\Lambda^m(\xi,\eta)=(\xi,\eta)+\Lambda^m(0^{l(\xi)},0^{l(\eta)})$. We define
\[ a_m(\xi,\eta)=\sum_{x,y \in
\Lambda^m(\xi,\eta)}\mbox{min}(x,y)-\sum_{x,y \in
\Lambda^m(0^{l(\xi)},0^{l(\eta)})}\mbox{min}(x,y).\] In this summation, we sum over
every pair of entries of the symbols. This function $a_m$ gives
rise to truncated induction: if $W'$ is a subgroup of $W_0$ and
$\chi'$ is a character of $W'$, let
\begin{equation}
\mbox{Ind}_{W'}^{W_0}(\chi')=\sum_{\chi \in
\hat{W_0}}n_{\chi',\chi} \chi.
\end{equation}
Suppose that $\chi_0 \in \hat{W}_0$ with $n_{\chi',\chi_0}>0$ is
such that $n_{\chi',\chi} >0$ implies $a_m(\chi)\leq a_m(\chi_0)$.
Then we define
\begin{equation}\label{trind}
\trm_{W'}^{W_0}(\chi')=\sum_{\chi:a_m(\chi)=a_m(\chi_0)}n_{\chi',\chi}\chi.
\end{equation}
In other words, we only keep the representations whose $a_m$-value
is maximal among the representations which occur in the induction.

\subsubsection{Definition of $\Sigma_m(W_0r)$} Let $W_0r$ be the
central character of an element of $\Hrcc$. We recall the
definition of Springer correspondents of $W_0r$, as in \cite{S}.
Let $q_2=q_1^m$. If $W_0r$ is a residual point, then we define its
Springer correspondents to be $\Sigma_m(W_0r)=[\S_m(\l)]_m$, where
$W_0r=W_0c(\l;k,mk))$. It is shown in \cite{S} that
$\Sigma_m(W_0r)$ is well-defined.

Otherwise, we can choose $r$ in its orbit $W_0r$ such that $r=r_L$
is the center of a real residual coset $L$ for which $R_L$ is a
standard parabolic root subsystem of $R_0$. Then $R_L$ is of type
$A_\kappa \times B_l:=A_{\kappa_1-1} \times A_{\kappa_2-1} \times
\dots A_{\kappa_r-1}\times B_l$ for some partition
$\kappa=(\kappa_1,\kappa_2,\dots,\kappa_r)$ of $n-l$. It is known
that an affine Hecke algebra of type $A_{p-1}$ and parameter $q$
has one Weyl group orbit of residual points (assuming that
$\mbox{rank}(X)=\mbox{rank}(Q)$), which is $W_0r_p$ with
\begin{equation}\label{r}
r_{p}(q)=(q^{-(p-1)/2},q^{-(p-3)/2},\dots,q^{(p-3)/2},q^{(p-1)/2}).
\end{equation}
The associated discrete series representation of $\H(A_{p-1},q)$
is one-dimensional, and its restriction to $\H_0(A_{p-1},q)$ is,
for $q>1$, the sign representation. Therefore, we associate to
$r_{p}$ a strip $S(p)$ with $\kappa$ boxes, which have entries
$(-(p-1)/2,-(p-3)/2,\dots,(p-3)/2,(p-1)/2)$ (the exponents of $q$
in \eqref{r}). We will call such a strip an $A$-strip. Notice that
its entries do not depend on $q$. We will also need the strip with
the absolute values of the same entries, we denote it by $|S(p)|$
and also call it a (positive) $A$-strip.

The residual point $W_Lr_L$ can be written as
\begin{equation}
\label{rL}
 r_L(q_1,q_2)=r_\kappa(q_1) c(\mu;q_1,q_2)=r_{\kappa_1}(q_1)r_{\kappa_2}(q_1)\dots r_{\kappa_r}(q_1)c(\mu;q_1,q_2),
\end{equation}
where
$\mu \vdash l$ is such that $c(\mu;q,q^m)$ is residual. By this notation, we mean that the first $\kappa_1$ coordinates of $r_L$ are as in \eqref{r} with $p=\kappa_1$, the next $\kappa_2$ coordinates are as in \eqref{r} with $p=\kappa_2$, etc., until finally the last $l$ coordinates are those of $c(\mu;q_1,q_2)$. Then we
define, letting $(\kappa_i)$ denote the trivial representation of $W_0(A_{\kappa_i-1})$: \[ \Sigma_m(W_Lr_L)=(\kappa_1) \otimes (\kappa_2) \otimes  \dots (\kappa_r)\otimes [\S_m(\mu)]_m=\mbox{triv}_\kappa \otimes [\S_m(\mu)]_m,\]
that is, we tensor the Springer correspondents of $W_0(B_l)c(\mu;q_1,q_2)$
with the Springer correspondents of $W_0(A_\kappa)r_\kappa(q_1)$.
Finally, we let
\begin{equation}
\label{defspringer} \Sigma_m(W_0r_L)=\{(\xi,\eta) \mid (\xi,\eta)
\mbox{ occurs in }\trm_{W_L}^{W_0}(\Sigma_m(W_Lr_L))\}.
\end{equation}
In other words, we truncatedly induce the Springer correspondents
of $W_Lr_L$ to obtain those of $W_0r_L$.

It is shown in \cite{S} that (under the identification of
bipartitions and irreducible characters of $W_0$) we have, for all
$m$, a disjoint union
\begin{equation}
\label{union}
 \bigcup_{L \mbox{ real residual coset}}\Sigma_m(W_0r_L) =
\hat{W}_0.
\end{equation}

It is not hard to show that all $\Sigma_m(W_0r_L)$ are singletons
unless $q_1,q_2$ are special parameters.

\subsection{Conjecture on $\Hrcc$} We can now formulate our conjecture on
$\Hrcc$. If $M$ is a $\H_0$-module, we denote by $R(M)$ the
corresponding $W_0$-module which we obtain in the limit
$q_{\ca}^{1/2} \to 1$.

\begin{conj}
\label{con} (\cite{S}) Let $\H$ be the affine Hecke algebra with
labels $q_2=q_1^m$ obtained from the $q_{\ca}$ as in \eqref{q1q2}
and $q_1>1,q_2 \geq 1$. Then we have the following description of
$\Hrcc$. Let $W_0r$ be the $W_0$-orbit of the center of a real
residual coset of $\H$.
\begin{itemize}
\item[(i)] $\Delta_{W_0r}$ is indexed by $\Sigma_m(W_0r)$.
\item[(ii)] The modules in $\Delta_{W_0r}$ are
naturally graded for the action of $\H_0$. The top non-zero degree of these modules is
$a_m(\Sigma_m(W_0r))$, and the $\H_0$-representations in this
degree are irreducible. The corresponding
$W_0$-representations are $\Sigma_m(W_0r) \otimes \e$.
\item[(iii)] Let $M_\chi \in \Hrcc$ be the $\H$-module labelled by $\chi \in \hat{W}_0$ using (ii) and \eqref{union}. That is, its top
degree is $R(M_\chi^{top})\cong \chi \otimes \e$. Let $\succ_m$ be
an ordering on $\hat{W_0}$ which is a refinement of the pre-order
given by $\psi \succ \chi \Rightarrow a_m(\psi) \leq a_m(\chi)$,
in which similarity classes form intervals. Then the matrix
$(R(M_{\chi}),\psi)_{\psi,\chi \in \hat{W}_0}$ is block-lower
triangular where the blocks are given by the similarity classes.
\item[(iv)] The $\H_0$-structure of $M_\chi$ can be computed with
the generalized Green functions of Shoji (see \cite{shoji}).
\end{itemize}
\end{conj}

The details of (iv), as well as the corresponding
 statements for
arbitrary $q_1,q_2$ can be found in \cite{S}. We remark that the conjecture is known to hold for $m=1$, by the theory of Kazhdan-Lusztig (cf. \cite{KL}).

\subsection{Reducibility of induced discrete series} \label{red}
In this article we will show that \ref{con}(i) for arbitrary
central character follows from \ref{con}(i) for residual points.
Indeed, Conjecture \ref{con} implies in particular the number of
irreducible components in an induced representation $\pi(\xi)=\pi(\Pi_L,\delta,1)$ for $\delta \in \hat{\H}^t_{L,\R}$
 as follows. Indeed, consider a residual coset $L$ such that $R_L$ has simple roots $\Pi_L$.
Let its center $r_L$ be as in \eqref{rL}. According to \ref{con}(i),
the cardinality of $\Delta_{W_Lr_L}$ is equal to the cardinality
of $\Sigma_m(W_Lr_L)$, and the cardinality of $\Delta_{W_0r_L}$ is
equal to the cardinality of $\Sigma_m(W_0r_L)$. It is shown in
\cite{S} that we have
\[ |\Sigma_m(W_0r_L)|=2^d |\Sigma_m(W_Lr_L)|,\]
where $d$ is the number of $A$-strips of length $k_i$
($i=1,2,\dots,r$) which can be glued to $T_m(\mu)$ while
maintaining the $m$-tableau of a partition. This number $d$ is
independent of $\mu'$ in the set $\{ \mu' \vdash l \mid \S_m(\mu')
\sim_m \S_m(\mu)\}$.

\begin{figure}[htb]
\begin{center}
{\includegraphics[angle=0,scale=0.3]{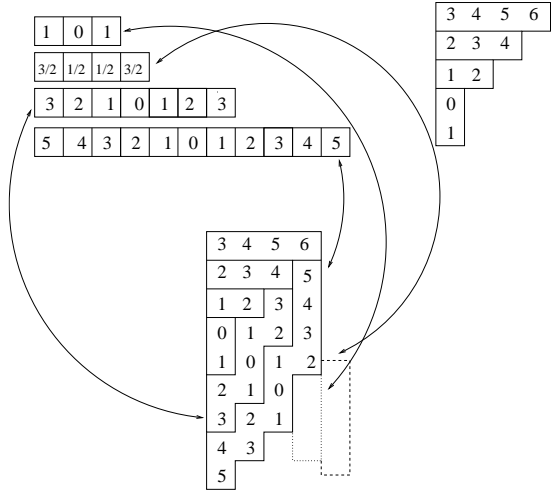}}
\end{center}
\caption{$m=3$, $\kappa=(3,4,7,11)$, $\mu=(1,1,2,3,4)$}
\label{example}
\end{figure}

As a consequence, we expect that $\pi(\xi)$ has
$2^d$ irreducible components. We give an example in Figure
\ref{example}, for $m=3$, $\kappa=(3,4,7,11)$ and
$\mu=(1,1,2,3,4)$. The upper part depicts the $A$-strips
$|S(3)|,|S(4)|, |S(7)|, |S(11)|$ and the $3$-tableau $T_3(\mu)$.
As shown in the lower part of the Figure, two of the strips (those
of $A_6$ and $A_{10}$) can be glued to $T_3(\mu)$, while this can
not be done for the strips of $A_2$ and $A_3$ (the dashed lines).
Therefore, $d=2$ and we expect $\pi(\xi)$ to have four irreducible
components.

\section{Explicit calculation of the $R$-group}

Let $R_L \subset R_0$ be a standard parabolic root system, then
$R_L$ has type $A_\kappa \times B_l=A_{\kappa_1-1} \times
A_{\kappa_2-1} \times \dots \times A_{\kappa_r-1} \times B_l$, for
some $l\geq 0$ and $\kappa=(\kappa_1,\dots,\kappa_r) \vdash n-l$.
The central character of a discrete series representation $\delta$
of $\H_L$ is then given by $W_Lr_L$, with $r_L$ as in \eqref{rL}.
Let $\kappa_0=0$, $\kappa_{r+1}=l$, and put
$K_i=\kappa_0+\kappa_1+\dots+\kappa_i$ for $i=0,1,\dots,r+1$. Define $L_i=\{\a_j \mid j=K_{i-1}+1,K_{i-1}+1,\dots,K_i-1\}$
($i=1,\dots,r$), $L_{r+1}=\{\a_j \mid j=K_r+1,\dots,n\}$ and
$I_i=\{K_{i-1}+1,\dots,K_i\}$. Then $I=\{1,\dots,n\}=\cup_i I_i$
and $L=\cup_iL_i$. Define $E_i=[\a_{K_{i-1}+1}]^L$ for
$i=1,2,\dots,r$, then $\{E_1,\dots,E_r\}$ forms a basis of $\t^L$.
We begin by calculating the root system $R_0(\xi)$.

\subsection{Calculation of $R_0(\xi)$}
One checks that we have the following relations between the projections of $R_0\backslash R_L$ onto $\t^L$:
\begin{itemize}
\item[(i-a)] Let $\a=e_{i_1} + e_{j_1}$ with $i_1 \in I_{p_1}, j_1 \in
I_{p_2}, p_1 \neq r+1, p_2 \neq r+1, p_1 \neq p_2$. Then
$[\a]^L=[e_{i_2}+ e_{j_2}]^L$ for all $i_2 \in I_{p_1}, j_2 \in
I_{p_2}$.
\item[(i-b)] Let $\a=e_{i_1} - e_{j_1}$ with $i_1 \in I_{p_1}, j_1 \in
I_{p_2}, p_1 \neq r+1, p_2 \neq r+1, p_1 \neq p_2$. Then
$[\a]^L=[e_{i_2}- e_{j_2}]^L$ for all $i_2 \in I_{p_1}, j_2 \in
I_{p_2}$.
\item[(ii)] Let $\a=e_{i_1}+e_{j_1}$ with $i_1,j_1 \in I_p$ with $p \neq
r+1$ (so that $\a \notin R_L$). Then
$[\a]^L=[e_{i_2}+e_{j_2}]^L=2[e_{i_2}]^L=2[e_{i_2} \pm e_{j_3}]^L$ for
all $i_2,j_2 \in I_p, j_3 \in I_{r+1}$. Clearly, $[e_{i_2}]^L$ is the primitive one.
\end{itemize}

To calculate the pole order of $c_\beta(r_L)$ for each primitive $\beta$, we note that
\[ e_{K_{i-1}+d}(r_L)=e_d(r_{\kappa_i})(q_1)=q_1^{-(\kappa_i-1)/2+(d-1)}, \ i=1,\dots,r;\ d=1,\dots,\kappa_i \]
and
\[ e_{K_r+i}(r_L)=e_{i}(c(\mu);q_1,q_2)=q_1^{c(\Box_i)}q_2; \ i=1,\dots,l.\]
for some numbering $\Box_1, \dots, \Box_l$ of the boxes of $Y(\mu)$, and $c(\Box)$ the content of $\Box$.

We are going to compute the pole order of $c_\beta(r_L)$ for every primitive $\beta$. Notice that, since the central character $W_Lr_L$ is real, this amounts to calculating the pole order of (see \eqref{cdef})
\begin{equation}
\label{cmodif}
C_{\beta}=\prod_{\{\a \mid [\a]^L \in \Z_{>0}\beta\}}\frac{1-q_{\ca/2}^{-1/2}q_{\ca}^{-1}\theta_{-\a}}{1-\theta_{-\a}},
\end{equation}
when evaluated in $r_L$. We calculate $C_\beta(r_L)$ for each case
above separately.

(i-a) Let $\a=e_i+e_j$ with $i \in I_{p_1}, j \in I_{p_2}$, $p_1 \neq p_2$, $p_1 \neq r+1$, $p_2 \neq r+1$. Then
\begin{equation}
\label{gevaltwee}
 C_\beta(r_L)=\prod_{d_1,d_2=1}^{\kappa_{p_1},\kappa_{p_2}}\frac{1-q_1^{-1}q_1^{-(\kappa_{p_1}-1)/2+(d_1-1)-(\kappa_{p_2}-1)/2+(d_2-1)}}{1-q_1^{-(\kappa_{p_1}-1)/2+(d_1-1)-(\kappa_{p_2}-1)/2+(d_2-1)}}
\end{equation}

We claim that \eqref{gevaltwee} has a pole only if
$\kappa_{p_1}=\kappa_{p_2}$. Indeed, \eqref{gevaltwee} has a pole
for every $(d_1,d_2)$ such that
$d_1+d_2=(\kappa_{p_1}+\kappa_{p_2})/2+1$ and a zero for every
$(d_1,d_2)$ such that $d_1+d_2=(\kappa_{p_1}+\kappa_{p_2})/2+2$.
Hence, if $\kappa_{p_1}+\kappa_{p_2}$ is odd then
\eqref{gevaltwee} does not have a pole. If
$\kappa_{p_1}+\kappa_{p_2}$ is even, then write
$\kappa_{p_2}=\kappa_{p_1}+2a$. The poles in \eqref{gevaltwee}
arise for
$\{(1,\kappa_{p_1}+a),(2,\kappa_{p_1}+a-1),\dots,(\kappa_{p_1},a+1)\}$
and the zeroes at
$\{(1,\kappa_{p_1}+a+1),(2,\kappa_{p_1}+a),\dots,(\kappa_{p_1},a+2)\}$,
except when $a=0$, when $(1,\kappa_{p_1}+1)$ does not occur in the
product. This proves the claim.

\medskip
(i-b) Suppose $\beta=[e_{i_1}-e_{j_2}]$ and $i_1 \in
I_{p_1}, j_1 \in I_{p_2}$, with $p_1 \neq p_2$, $p_1 \neq r+1$, $p_2 \neq r+1$. Then
\[ C_\beta(r_L)=\prod_{d_1,d_2=1}^{p_1,p_2}\frac{1-q_1^{-1}q_1^{(\kappa_{p_1}-\kappa_{p_2})/2-(d_1-d_2)}}{1-q_1^{(\kappa_{p_1}-\kappa_{p_2})/2-(d_1-d_2)}}\]
We claim that this expression has a pole only for
$\kappa_{p_1}=\kappa_{p_2}$. The proof is identical to the one
above: if $\kappa_{p_1}+\kappa_{p_2}$ is odd then there are no
poles at all, hence suppose that $\kappa_{p_2}\geq \kappa_{p_1}$,
say $\kappa_{p_2}=\kappa_{p_1}+2a$ for some $a \in \Z_{\geq 0}$.
Then the poles are indexed by $(d_1,d_2) \in
\{(a+1,1),(a+2,2),\dots,(a+\kappa_{p_1},\kappa_{p_1})\}$ and the
zeroes by $(d_1,d_2) \in
\{(a,1),(a+1,1),\dots,(a+\kappa_{p_1}-1,\kappa_{p_1})\}$, except
when $a=0$ (i.e., when $\kappa_{p_1}=\kappa_{p_2}$), in which case
$(a,1)$ does not occur.

\medskip
(ii) Let $\a=e_i$ for $i \in I_p$, $p \neq r+1$ and let $\beta=[\a]^L$. Then
\begin{equation}
\label{c}
C_\beta(r_L)=\prod_{d=1}^{\kappa_p}\frac{1-q_2^{-1}q_1^{-(\kappa_p-1)/2+(d-1)}}{1-q_1^{-(\kappa_p-1)/2+(d-1)}}
\prod_{1\leq d_1<d_2 \leq
\kappa_p}\frac{1-q_1^{-1}q_1^{-\kappa_p+d_1+d_2-1}}{1-q_1^{-\kappa_p+d_1+d_2-1}}\times
\end{equation}
\[ \times \prod_{d=1}^{\kappa_p}
\prod_{i=1}^l\frac{1-q_1^{-1}q_1^{-(\kappa_p-1)/2+(d-1)}q_1^{c(\Box_i)}q_2}{1-q_1^{-(\kappa_p-1)/2+(d-1)}q_1^{c(\Box_i)}q_2}\frac{1-q_1^{-1}q_1^{-(\kappa_p-1)/2+(d-1)}q_1^{-c(\Box_i)}q_2^{-1}}{1-q_1^{-(\kappa_p-1)/2+(d-1)}q_1^{-c(\Box_i)}q_2^{-1}}.\]
We write this expression as
\begin{equation}
\label{c2}
C_\beta(r_L)=C_\beta(\kappa_p) C_\beta(\kappa_p,T_m(\mu)),
\end{equation}
where $C_\beta(\kappa_p)$ denotes the top line of \eqref{c} and
$C_\beta(\kappa_p,T_m(\mu))$ the bottom line. Thus,
$C_\beta(\kappa_p)$ denotes the part of \eqref{c} which
corresponds to the factor $A_{\kappa_p-1}$, and
$C_\beta(\kappa_p,T_m(\mu))$ denotes the part of \eqref{c} which
calculates the ``interaction'' between $A_{\kappa_p-1}$ and
$T_m(\mu)$.

\medskip
(ii-a) First we calculate the pole order of $C_\beta(\kappa_p)$. In the left hand product of $C_\beta(\kappa_p)$, we find a pole if and only if $\kappa_p$ is odd. In the right hand side
product, we have a pole for $(d_1,d_2)\in
\{(1,\kappa_p),(2,\kappa_p-1),\dots,(\lfloor\frac{\kappa_p}{2}\rfloor,\lceil\frac{\kappa_p}{2}\rceil+1)\}$.
The total number of poles in $C_\beta(\kappa_p)$ is therefore $\lceil \frac{\kappa_p}{2}\rceil$.

Next, we count the number of zeroes in $C_\beta(\kappa_p)$. In the left hand product, we find one zero factor if $q_2=q_1^m$ for  $m \in \{ (\kappa_p-1)/2,(\kappa_p-3)/2,\dots\}$, and none otherwise. In the right
hand side we have, independently of the relation between $q_1$ and $q_2$, a zero for every $(d_1,d_2) \in
\{(2,\kappa_p),(3,\kappa_p-1),\dots,(\lfloor\frac{\kappa_p+1}{2}\rfloor,\lceil\frac{\kappa_p+3}{2}\rceil)\}$, i.e., a total on the right hand side of $\lfloor\frac{\kappa_p+1}{2}\rfloor-1=\lceil\frac{\kappa_p}{2}\rceil-1$ zeroes.

Thus,
\begin{equation}
\label{c2boven}
 C_\beta(\kappa_p){\rm\ has\ a\ pole\ of \ order\ } \left\{ \begin{array}{ll} 0 & {\rm if}\ q_2=q_1^m, m \in \{ (\kappa_p-1)/2,(\kappa_p-3)/2,\dots\},\\ 1 & {\rm  otherwise}.\end{array}\right.
\end{equation}
%In particular, \marginpar{deze zin weg?}if $q_2=1$ then
%$C_\beta(\kappa_p)$ has a pole if $\kappa_p$ is even and does not
%have a pole if $\kappa_p$ is odd.
\medskip

(ii-b) Now we calculate the pole order of
$C_\beta(\kappa_p,T_m(\mu))$. Suppose that $q_2=q_1^m$ with $m \in
\half\Z$ (otherwise the parameters are generic and
$C_\beta(\kappa_p,T_m(\mu))$ does not have a pole). Then
$C_\beta(\kappa_p,T_m(\mu))=$
\begin{equation}
\label{onder} \prod_{d=1}^{\kappa_p}
\prod_{i=1}^l\frac{1-q_1^{-1}q_1^{-(\kappa_p-1)/2+(d-1)}q_1^{c(\Box_i)+m}}{1-q_1^{-(\kappa_p-1)/2+(d-1)}q_1^{c(\Box_i)+m}}\frac{1-q_1^{-1}q_1^{-(\kappa_p-1)/2+(d-1)}q_1^{-(c(\Box_i)+m)}}{1-q_1^{-(\kappa_p-1)/2+(d-1)}q_1^{-(c(\Box_i)+m)}}.
\end{equation}
Notice that the $|c(\Box)+m|$ are the entries of the $m$-tableau
$T_m(\mu)$.

%Graphically, \eqref{c} having or not having a pole should correspond to whether or not it is possible to glue an $A$-strip of length $\kappa_p$ to $T_m(\mu)$.

Let $S(\kappa_p)$ be the $A$-strip with
entries $(-(\kappa_p-1)/2,-(\kappa_p-3)/2,\dots,(\kappa_p-3)/2,(\kappa_p-1)/2)=(-z,\dots,z)$. For $\Box \in S(\kappa_p)$ a square of the $A$-strip, let $e(\Box)$ be its entry. So, $e(\Box_i)=e_{K_{p-1}+i}(r_L)=q_1^{-z+i-1} (i=1,\dots,\kappa_p)$. For a square $\Box$ of $Y(\l)$, let $e_m(\Box)$ be its entry in $T_m(\l)$. With these notations,
\begin{equation}
\label{onder2}
C_\beta(\kappa_p,T_m(\mu))=\prod_{\Box \in S(\kappa_p)}\prod_{\Box' \in T_m(\mu)} \frac{1-q_1^{-1}q_1^{-e(\Box)+e_m(\Box')}}{1-q_1^{-e(\Box)+e_m(\Box')}} \frac{1-q_1^{-1}q_1^{-e(\Box)-e_m(\Box')}}{1-q_1^{-e(\Box)-e_m(\Box')}}.
\end{equation}
Thus, $C_\beta(\kappa_p,T_m(\mu))$ has a pole for every $\Box \in
S(\kappa_p)$, $\Box'\in T_m(\mu)$ such that $e(\Box)=\pm
e_m(\Box')$ (where, if $e_m(\Box')=0$, we count this equality
twice); and a zero for every $\Box \in S(\kappa_p), \Box' \in
T_m(\mu)$ such that $e(\Box)\pm e_m(\Box')=1$ (where again, we
count twice if $e_m(\Box')=0$). Now recall the splitting map
$\S_m$ which divides $T_m(\mu)$ into horizontal and vertical
blocks. We write $B \in \S_m(\mu)$ to denote that $B$ is one of
the blocks into which $T_m(\mu)$ is partitioned by $\S_m$. We can
thus write
\begin{eqnarray*}
\label{onder3}
C_\beta(\kappa_p,T_m(\mu))&=&\prod_{\Box \in S(\kappa_p)}\prod_{B \in \S_m(\mu)}\prod_{\Box' \in B} \frac{1-q_1^{-1}q_1^{-e(\Box)+e_m(\Box')}}{1-q_1^{-e(\Box)+e_m(\Box')}} \frac{1-q_1^{-1}q_1^{-e(\Box)-e_m(\Box')}}{1-q_1^{-e(\Box)-e_m(\Box')}}\\&=& \prod_{B \in \S_m(\mu)}C_\beta(\kappa_p,B).
\end{eqnarray*}
We can thus calculate the pole-order of $C_\beta(r_L)$ as the sum of the pole orders of $C_\beta(\kappa_p)$ and those of the $C_\beta(\kappa_p,B)$. We recall from \cite{S} that we may choose $\mu$ such that the splitting procedure $\S_m$ first selects $m$ horizontal blocks, and then alternatingly a vertical and a horizontal block. This means that we can suppose that a block $B$ contains entries $(x,x+1,\dots,y)$ with $0\leq x \leq y$. However, $x=y=0$ does not occur since $T_m(\mu)$ does not contain a block of length one containing a zero.

Denote by $B$ one of these blocks.
%First we consider the case where $B$ is one of the first $m$ blocks emerging from the splitting procedure, i.e., $B$ contains a sequence $(x,x+1,\dots,y)$ where $y \geq x$ and $1\leq x \leq m$.
 Then the pole order of $C_\beta(\kappa_p,B)$ is given by $|P(\kappa_p,B)|-|Z(\kappa_p,B)|$, where
\[ P(\kappa_p,B)=\{(e(\Box),e_m(\Box'),\varepsilon) \mid \Box \in S(\kappa_p), \Box' \in B, \e \in \{+,-\}, e(\Box)+\varepsilon e_m(\Box')=0\},\]
\[ Z(\kappa_p,B)=\{(e(\Box),e_m(\Box'),\varepsilon) \mid \Box \in S(\kappa_p), \Box' \in B, \e \in \{+,-\}, e(\Box)+\varepsilon e_m(\Box')=1\}.\]

%We will denote poles and zeroes of \eqref{onder} by a pair $(e(\Box),e_m(\Box'))$ where $\Box \in S(\kappa_p)$ and $\Box' \in B \subset T_m(\mu)$. The set of zeroes is denoted by $Z(B,S(\kappa_p))=\{ (e(\Box),e_m(\Box')) \mid \Box \in S(\kappa_p), \Box' \in B, e(\Box)\pm e_m(\Box')=1\}$; likewise the set of poles is denoted $P(B,S(\kappa_p))$.
We calculate $P(\kappa_p,B)$ and $Z(\kappa_p,B)$ by a case by case analysis. Recall $z=(\kappa_p-1)/2$. Clearly, if $x-z \notin \Z$ then $P(\kappa_p,B)=Z(\kappa_p,B)=\emptyset$ since all $e(\Box)\pm e_m(\Box') \notin \Z$. Otherwise,
\begin{itemize}
\item If $z<x-1$ then $P(\kappa_p,B)=Z(\kappa_p,B)=\emptyset$.
\item If $z=x-1$ then
$Z(\kappa_p,B)=\{(-z,x,+)\}$ and $P(\kappa_p,B)=\emptyset$.
\item If $x\leq z<y$ then $P(\kappa_p,B)=\{(\pm z,z,\mp),(\pm(z-1),z-1,\mp),\dots,(\pm x,x,\mp)\}$ and $Z(\kappa_p,B)=\{(-z,z+1,+),(-z+1,z,+),\dots,(-x+1,x,+)\}\cup \{(z,z-1,-),(z-1,z-2,-),\dots,(x+1,x,-)\}$. One checks that $|P(\kappa_p,B)|=|Z(\kappa_p,B)|=2(z-x+1)$.
\item If $x\leq y=z$, then $P(\kappa_p,B)=\{(\pm z,z,\mp),(\pm (z-1),z-1,\mp),\dots,(\pm x,x,\mp)\}$ and $Z(\kappa_p,B)=\{(-z+1,z,+),(-z+2,z-1,+),\dots,(-x+1,x,+)\}\cup\{(z,z-1,-),(z-1,z-2,-),\dots,(x+1,x,-)\}$.
Thus, we find a pole of order one.
\item If $z>y$, then $P(\kappa_p,B)=\{(\pm y,y,\mp),(\pm (y-1),y-1,\mp),\dots,(\pm x,x,\mp)\}$ and $Z(\kappa_p,B)=\{(-y+1,y,+),(-y+2,y-1,+),\dots,(-x+1,x,+)\}\cup\{(y+1,y,-),(y,y-1,-),\dots,(x+1,x,-)\}$. In this case, $|P(\kappa_p,B)|=|Z(\kappa_p,B)|$.
\end{itemize} Summarizing, we have:

\begin{lemma}\label{plak}
Let $z=(\kappa_p-1)/2$. Let
$B \in \S_m(\mu)$ be a block with entries $(x,x+1\dots,y)$ with $0
\leq x \leq y$. Then $C_\beta(\kappa_p,B)$ has pole order
\[ \begin{array}{c@{\mbox{ if }}c} -1 & z=x-1 \\ 0 & z \notin \{x-1,y\} \\ 1 &
z=y. \end{array}\]
\end{lemma}

This enables us to show

\begin{prop} Let $(\Pi_L,\delta,1) \in \Xi$ have real central character $W_Lr_L$ with $r_L$ as in \eqref{rL}. Consider the primitive vector $\beta=E_p=[e_{K_{p-1}+1}]^L \in \t^{L^*}$. %corresponding to the factor $A_{\kappa_p-1}$ in the standard parabolic $R_L\subset R_0$.
Then $C_\beta$ has a pole in $r_L$ unless a (positive) $A$-strip $|S(\kappa_p)|$ of length $\kappa_p$ can be glued to $T_m(\mu)$.
\end{prop}

\pf: Let $z=(\kappa_p-1)/2$. If $m-z \notin \Z$, then $S(\kappa_p)$ can certainly not be glued to $T_m(\l)$. Recall that we decompose $C_\beta(r_L)=C_\beta(\kappa_p)C_\beta(\kappa_p,T_m(\mu))$. By \eqref{c2boven}, $C_\beta(\kappa_p)$ has a pole whereas by Lemma \ref{plak}, $C_\beta(\kappa_p,T_m(\mu))$ has pole order zero. Thus $C_\beta(r_L)$ has a pole of order one as required.

We therefore assume that $m-z \in \Z$. If $l=0$, then $T_m(\mu)$ is the empty tableau, so $C_\beta(r_L)=C_\beta(\kappa_p)$. One can place a strip of length $\kappa_p$ on the empty tableau if and only if $(\kappa_p-1)/2 \geq m$. By \eqref{c2boven}, this can be done if and only if $C_\beta(\kappa_p)$ does not have a pole.

Now suppose $l>0$.

(i) Suppose that $|S(\kappa_p)|$ can not be glued to $T_m(\l)$. Then there are two possibilities.

(i-a) We have $z<m-l(\mu)$. In particular, $l(\mu)<m$ and thus $T_m(\mu)$ contains only horizontal blocks: we have $\S_m(\mu)=(\mu,-)$. These blocks have initial entries $m,m-1,\dots,m-l(\mu)+1$. Since $z<m$, $C_\beta(\kappa_p)$ has a pole by \eqref{c2boven}. By Lemma \ref{plak}, $C_\beta(\kappa_p,T_m(\mu))$ has pole order zero, so $C_\beta$ has a pole of order one in $r_L$.

(i-b) $T_m(\mu)$ contains a (unique) block $B$ whose last entry is $z$. In this case, $C_\beta(\kappa_p,B)$ has a pole of order one. If $m>z$ then $C_\beta(\kappa_p)$ has a pole, but $m>z$ implies that there exists a block $B'$ of $T_m(\mu)$ which starts on $z+1$. If $z>0$ then this block is unique and the total pole order of $C_\beta(r)$ is one. If $z=0$, then $B$ has last entry zero, which is impossible. If $m \leq z$ then $C_\beta(\kappa_p)$ does not have a pole. But $m \leq z$ implies that no block starts on $z+1$, so again $C_\beta(r)$ has a pole of order one.

(ii) Now suppose that $|S(\kappa_p)|$ can be glued to $T_m(\l)$. Then no block in $T_m(\mu)$ ends on $z$. If $m \leq z$, then by \eqref{c2boven}, $C_\beta(\kappa_p)$ has pole order zero. There is no block $B'$ in $T_m(\mu)$ which starts on $z+1 \geq m+1$, hence $C_\beta(\kappa_p,B')$ has pole order zero for all $B'$. Thus, $C_\beta(r_L)$ does not have a pole.

%If $m=z$ then $C_\beta(\kappa_p)$ does not have a pole. The tableau $T_m(\mu)$ does not contain a block which starts on $z+1>m$ so the total pole order remains zero.

If $m > z$, then $C_\beta(\kappa_p)$ has a pole. But $m > z$ implies that $T_m(\mu)$ contains a block $B'$ which starts on $z+1 \leq m$, and so $C_\beta(\kappa_p,B')$ cancels the pole of $C_\beta(\kappa_p)$. Hence, $C_\beta(r_L)$ does not have a pole. \qed

\subsubsection{Description of $R_0(\xi)$}
We can now describe the root system $R_0(\xi)$. We have seen that $E_i
\pm E_j \in R_0(\xi)$ if and only if $\kappa_i=\kappa_j$, and $E_i
\in R_0(\xi)$ unless an $A$-strip $|S(\kappa_i)|$ can be glued to
$T_m(\mu)$. Therefore, we have

\begin{prop} Write \begin{equation}
\label{hergroep}
(\kappa_1,\dots,\kappa_r)=(l_1^{r_1},\dots,l_s^{r_s})
\end{equation}
by grouping the equal $\kappa_i$ together. Then the root system $R_0(\xi)$ is of type
\begin{equation}
\label{algemeen}
 R_0(\xi) \simeq R_0(X_{r_1}) \times R_0(X_{r_2}) \times
\dots \times R_0(X_{r_s}),
\end{equation}
 where
\[ X_{r_i}=\left\{ \begin{array}{c@{\mbox{ if }}l} D_{r_i} & \mbox{a strip of
length } l_i \mbox{ fits into }T_m(\mu),\\ B_{r_i} &
\mbox{not}.\end{array}\right.\]
\end{prop}
In this formula, we use the convention that $R_0(D_1)$ is the empty root system, while for $n \geq 2$, the root system of type $D_n$ on an orthogonal basis $\{e_1,\dots,e_n\}$ is given by the vectors $\{\pm e_i \pm e_j\mid i \neq j \}$.

\subsection{Calculation of $R(\xi)$}
We are going to compute
\[ R(\xi)=\{ w \in W_0 \mid w(L)=L, w(R_0^+(\xi))=R_0^+(\xi)\}.\]
Observe that since $w(L)=L$, it follows that if $w \in R(\xi)$ and
if $w(L_{i_1})=L_{i_2}$ (which is always true for some $i_2$), then
$\kappa_{i_1}=\kappa_{i_2}$. Thus, $w \in R(\xi)$ automatically satisfies
$w(R_0^+(X_{r_j}))=R_0^+(X_{r_k})$, for $j,k$ such that
$r_j=r_k$. We show that moreover, $j=k$.

\begin{lemma}\label{vast}
Let $w \in R(\xi)$ be such that $w(R_0^+(X_{r_j}))=R_0^+(X_{r_k})$ where $X\in \{B,D\}$. Then $j=k$.
\end{lemma}

\pf: Suppose $j \neq k$. By construction, $l_j \neq l_k$. Let
$r_j=r_k=p$. Then we may suppose that the root systems $R_0(X_p)$ are realised on
basis vectors $E_1,\dots,E_p$ (for $\kappa_1=\dots=\kappa_p=l_j$) and $E_{p+1},\dots,E_{2p}$ (for $\kappa_{p+1}=\dots=\kappa_{2p}=l_k$).

If $X=D$, then $w(\{E_i \pm E_j \mid 1 \leq i < j \leq p\})=\{E_i \pm E_j \mid p+1 \leq i< j \leq 2p\}$. Moreover, $w \in W_0$, so $w(E_i) \in \{E_1,\dots,E_{2p}\}$ for all $1 \leq i \leq 2p$. This means that we have $w(E_i)=E_{p+i}$ for all $1 \leq i <p$ and $w(E_p)=\pm E_{2p}$.

Thus, $[w(e_1)]^L=E_{p+1}$, so $w(e_1)=e_t$ for $t \geq pl_j+1$. But $w(e_1)=w(\a_1)+w(e_2)$, and $w(\a_1)=\a_s$ with $s<pl_j$. Thus, $w(e_2)=e_t-\a_s \notin R_0$, which is a contradiction.

If $X=B$, then similarly, we find that $w(E_i)=E_{p+i}$ for all $1 \leq i \leq p$ which also leads to a contradiction. Thus, $j=k$.
\qed

It follows that if $w \in R(\xi)$, then $w(L_j)=L_j$ for all $j$ (i.e., the parts of equal length in $L$ are not interchanged). In particular, if $w \in R(\xi)$ then $w(E_i)=\pm E_i$ for all $i=1,2,\dots,r$.

\begin{prop}
\label{Rgen}
If $R_0(\xi)$ contains only root systems of type $B$, then $R(\xi)=\{1\}$.
\end{prop}

Let $w \in R(\xi)$. We want to show that $w=1$, and therefore it
suffices to show that $w=1$ on all $I_i$ separately, which we
prove using (downward) induction on $i$.

(i) Suppose $i=r+1$. Then $R_{L_{r+1}}$ is a root system of type $B_l$, hence $w(L_{r+1})=L_{r+1}$ implies $w=1$ on $R_{L_{r+1}}$.

%(ii) Suppose $i=r$. Then $L_r=\{\a_{K_{r-1}+1},\dots,\a_{n-l-1}\}$ (possibly empty) and $w(L_r)=L_r$. Thus, $w(e_{K_{r-1}+1})=w(\a_{K_{r-1}+1}+\dots+\a_n)=\a_{K_{r-1}+1}+\dots+\a_{n-l-1}+w(\a_{n-l})+\a_{n-l+1}+\dots+\a_n$.
%On the other hand, since $[e_{K_{r-1}+1}]^L=E_r$ and $w(E_r)=E_r$, we have $w(e_{K_{r-1}+1})=\a_i+\dots+\a_n$ for some $K_{r-1}+1 \leq i \leq K_r=n-l$. Since $w(e_{K_{r-1}+1}) \in R_0$, it follows that $w(\a_{n-l})=\a_{n-l}$. But then $w(L_r)=L_r$ and $w(\a_{n-l})=\a_{n-l}$ imply that $w(L_r)=L_r$ pointwise. Thus, we have $w=1$ on $I_r$.

(ii) Now suppose that $w=1$ on $\cup_{j>i+1}I_j$. Then we consider
$e_{K_i+1}=\a_{K_i+1}+\dots+\a_n$. One has
$w(e_{K_i+1})=w(\a_{K_i+1}+\a_{K_{i}+2}+\dots+\a_{K_{i+1}-1}+\a_{K_{i+1}}+\a_{K_{i+1}+1}+\dots+\a_n)=\a_{K_i+1}+\a_{K_{i}+2}+\dots+\a_{K_{i+1}-1}+w(\a_{K_{i+1}})+\a_{K_{i+1}+1}+\dots+\a_n$.
On the other hand, since $w(E_i)=E_i$, we have $w(e_{K_i+1})=e_j$ for some $j \in \{K_i+1,\dots,K_{i+1}\}$. Now,
$w(\a_{K_{i+1}})=w(e_{K_i+1})-(\a_{K_i+1}+\dots+\a_{K_{i+1}-1})-(\a_{K_{i+1}+1}+\dots+\a_n)=e_j-(e_{K_i+1}-e_{K_{i+1}})-e_{K_{i+1}+1}$
must be a root. Therefore, $j=K_i+1$, $w(e_{K_i+1})=e_{K_i+1}$ and $w(\a_{K_{i+1}})=\a_{K_{i+1}}$.
Since we already knew that $w(L_{i+1})=L_{i+1}$, it follows that $w(I_{i+1})=I_{i+1}$. %w(L_{i+1} \cup \{\a_{K_{i+1}}\})=w(\{\a_{K_i+1},\a_{K_i+2},\dots,\a_{K_{i+1}}\})=\{\a_{K_i+1},\a_{K_i+2},\dots,\a_{K_{i+1}}\}$
 But then $w(\a_{K_{i+1}})=\a_{K_{i+1}}$ implies that $w=1$ on $I_{i+1}$.
This proves the induction step and therefore the claim.\qed

It remains to calculate $R(\xi)$ for $R_0(\xi)$ which contains a root system of type $D$. The idea is that $w$ in $R(\xi)$ can interchange the
order of the simple roots on the last part of $L$ which gives a
basis vector for the type $D$ factor in \eqref{algemeen}. For
example, if $L$ is
\[ \bullet-\bullet-\circ-\bullet-\bullet-\circ-\bullet-\bullet
=\circ\] then $\t^L$ has basis vectors $E_1,E_2,E_3$. If $q_2=q_1$, then $R_0^+(\xi)=\{E_1\pm E_2, E_1 \pm E_3, E_2 \pm E_3\}$ (cf. \eqref{algemeen}). As seen in the proof of Lemma \ref{vast}, an element $w \in R(\xi)$ must fix $E_1$ and $E_2$ and may send $E_3$ to $\pm E_3$. Together with $w(L)=L$, we will see that this implies that if $w(E_3)=E_3$ then $w=1$, and if $w(E_3)=-E_3$ then $w$ interchanges the last two simple roots in $L$ ($\a_7$ and $\a_8$).

\begin{thm}\label{prop1} Let $\xi=(\Pi_L,\delta,1) \in \Xi$ have real central character $W_Lr_L$ as in \eqref{rL}. Suppose $R_L$ is of type $A_\kappa \times B_l$. Then \[ R(\xi) \simeq \Z_2^d,\] where $d$ is the number of $A$-strips $|S(\kappa_i)|$ which can be glued to $T_m(\l)$.
\end{thm}

\pf: By conjugating, we may assume that the root systems of type
$D$ are in the beginning of the Dynkin diagram (i.e., arise on
$\kappa_1,\kappa_2,\dots$), so by Proposition \ref{Rgen} we may assume that
$w(\a_i)=\a_i$ for all $\a_i$ not in these $D$-factors. Otherwise
stated, we can suppose for the calculation that $R_0(\xi)$ only
contains $D$-factors.

By Lemma \ref{vast}, if $w \in R(\xi)$, then $w$ fixes the positive roots in each irreducible factor of $R_0(\xi)$ separately.

(1) We first consider the case where all $A$-factors are of equal
rank, i.e., $\kappa=(l_1^{r_1})=(l^r)$ and
$R_0(\xi)=R_0(D_{r})$. Recall the basis vectors $E_1,\dots,E_r$ of $\t^L$. Suppose $w \in R(\xi)$. Then, as seen in (the proof of) Lemma \ref{vast}, we have % $w(\{E_i\pm E_j\mid 1\leq i<j\leq r\})=\{ E_i \pm E_j\mid 1\leq i <j\leq r \}$. Since the $E_i$ are determined as linear
%combinations of the $E_i\pm E_j$ and v.v., it is equivalent to
%calculate $w$ on the $E_i$. Now $w$ does not change the lengths,
%so one must have $w(E_i)=E_j$ for some $j$ and this implies that
 $w(E_i)=E_i$ for all $i<r$ and $w(E_r)=\pm E_r$. If $w(E_r)=E_r$, we have seen in Proposition \ref{Rgen} that $w=1$. We thus assume that $w(E_r)=-E_r$ and show that this uniquely determines $ 1\neq w\in R(\xi)$. We proceed in several steps.

(1a) Let $\beta_r=e_{(r-1)l+1}$, that is, $\beta_r$ is the first
coordinate vector of the last copy of $A_{l-1}$ in $R_L$. Then
$[\beta_r]^L=E_r$, hence $[w(\beta_r)]^L=-[\beta_r]^L$. This means
that $w(\beta_r)=-e_i$ for some $(r-1)l+1\leq i \leq r l$. We
have, since $w(L_i)=L_i$ for all $L_i$:
\begin{eqnarray*}w(\beta_r)&=&w(\a_{(r-1)l+1}+\dots+\a_{r l
-1}+\a_{r l}+e_{r l+1})
\\&=&\a_{(r-1)l+1}+\dots+\a_{r l -1}+w(\a_{r l})+e_{r
l+1}\\&=& \beta_r+w(\a_{r l})-\a_{r l}.\end{eqnarray*} Thus, $R_0
\ni w(\a_{r l})=w(\beta_r)-\beta_r+\a_{r l}=-e_i-e_{(r-1)l+1}+e_{r
l}-e_{r l +1}$, which implies that $i=r l$. Thus,
$w(\beta_r)=-e_{r l}$ and $w(\a_{r l})=-e_{(r-1)l+1}-e_{r l +1}$.
This implies the values of $w$ on $L_r$, since $w(\a_{r l-1}+\a_{r
l})=w(\a_{r l-1})+w(\a_{r l})=e_i-e_{i+1}-e_{(r-1)l+1}-e_{r l+1}$
for $i \in \{(r-1)l+1,r l-1\}$. It follows that $w(\a_{r
l-1})=\a_{(r-1)l+1}$, and analogously that $w(\a_{r
l-j})=\a_{(r-1)l+j}$ for all $j$. Thus, $w$ acts on $L_r$ by
flipping it to its mirror image.

%If there are only factors of type $A$ in $L$ then $r \kappa=n$, otherwise $r \kappa<n$ and $w$ fixes the $B$-factor in $L$ pointwise.

(1b) Let $\beta_{r-1}=e_{(r-2)l+1}$, that is, $\beta_{r-1}$ is the
first coordinate vector of the last-but-one copy of $A_{l-1}$ in
$R_L$. Then $[\beta_{r-1}]^L=E_{r-1}$, hence
$[w(\beta_{r-1})]^L=[\beta_{r-1}]^L$. This means that
$w(\beta_{r-1})=e_i$ for some $(r-2)l+1\leq i \leq (r-1) l$. We
have, since $w(L_i)=L_i$ for all $L_i$:
\begin{eqnarray*}w(\beta_{r-1})&=&w(\a_{(r-2)l+1}+\dots+\a_{(r-1) l
-1}+\a_{(r-1) l}+\beta_r)
\\&=&\a_{(r-2)l+1}+\dots+\a_{(r-1) l -1}+w(\a_{(r-1) l})-e_{r l}\end{eqnarray*} Thus, $R_0
\ni w(\a_{(r-1) l})=w(\beta_{r-1})-e_{(r-2)l+1}+e_{(r-1) l}+e_{r
l}$, which implies that $i=(r-2) l+1$. Thus,
$w(\beta_{r-1})=\beta_{r-1}$ and $w(\a_{(r-1) l})=e_{(r-1)l}+e_{r
l}$. This implies the values of $w$ on $L_{r-1}$, since
$w(\a_{(r-1) l-1}+\a_{(r-1) l})=w(\a_{(r-1) l-1})+w(\a_{(r-1)
l})=e_i-e_{i+1}+e_{(r-1)l}+e_{r l}$ for $i \in \{(r-2)l+1,(r-1)
l-1\}$. It follows that $w(\a_{(r-1) l-1})=\a_{(r-1)l-1}$, and
analogously that $w(\a_{(r-1) l-j})=\a_{(r-1)l-j}$ for all $j$.
Thus, $w$ acts on $L_{r-1}$ as the identity.

(1c)  Let $\beta_{r-2}=e_{(r-3)l+1}$. Then
$[\beta_{r-2}]^L=E_{r-2}$, hence
$[w(\beta_{r-2})]^L=[\beta_{r-2}]^L$. This means that
$w(\beta_{r-2})=e_i$ for some $(r-3)l+1\leq i \leq (r-2) l$. We
have, since $w(L_i)=L_i$ for all $L_i$:
\begin{eqnarray*}w(\beta_{r-2})&=&w(\a_{(r-3)l+1}+\dots+\a_{(r-2) l
-1}+\a_{(r-2) l}+\beta_{r-1})
\\&=&\a_{(r-3)l+1}+\dots+\a_{(r-2) l -1}+w(\a_{(r-2) l})+e_{(r-2)l+1}\end{eqnarray*} Thus, $R_0
\ni w(\a_{(r-2) l})=w(\beta_{r-2})-e_{(r-3)l+1}+e_{(r-2)
l}-e_{(r-2) l+1}$, which implies that $i=(r-3) l+1$. Thus,
$w(\a_{(r-2) l})=e_{(r-2)l}+e_{(r-2) l+1}$ and
$w(\beta_{r-2})=\beta_{r-2}$. This implies the values of $w$ on
$L_{r-2}$, as before. We find that $w$ acts on $L_{r-1}$ as the
identity.

(1d) It is clear that by repetition of this argument, one finds
that $w$ is uniquely determined, and that $w(L_i)=L_i$ pointwise
on every $i=1,\dots,r-1$. From the description of $w$ in (1a-c)
one sees easily that $w(e_i)=e_i$ except when $(r-1)l+1 \leq i
\leq rl$, for which we have $w(e_{(r-1)l+j})=-e_{rl-(j-1)}$
($j=1,\dots,l$). It follows that $w$ is an involution, and so
$R(\xi) \simeq \Z_2$.

(2) Now we consider the general case where $\kappa=(l_1^{r_1}\dots
l_s^{r_s})$. Let $N_i=\sum_{j=1}^{i-1} l_is_i$ for $2 \leq i \leq
s+1$ and $N_1=0$. Then $e_{N_i+1}$ corresponds to the first
coordinate vector of the first factor of type $A_{l_i-1}$. In part
(1) we have seen that $w \in R(\xi)$ satisfies
$w(e_{N_s+1})=e_{N_s+1}$ (if $r_s>1$) or
$w(e_{N_s+1})=-e_{N_{s+1}}$ (if $r_s=1$).  It is easy to check
that this %$w(e_1)=e_1$ (if $r \geq 2$) or $w(e_1)=-e_l$ (if $r=1$)
implies that we can repeat the arguments (1a-d) on each root
system of type $D$ in \eqref{algemeen} separately, to obtain for
every root system of type $D_{r_i}$ in $R_0(\xi)$ (including the
empty ones, which correspond to $r_i=1$) an involution $w_{i}$. If
$l_i>2$ then $w_{i}$ behaves as the non-trivial diagram
automorphism on the last copy of $A_{l_i}$ and as the identity on
the others. It is easy to see from the construction of the $w_{i}$
that they commute. We have thus found, for every $l_i$ such that a
strip $|S(l_i)|$ can be glued to $T_m(\mu)$, an involution $w_i
\in
R(\xi)$ and these commute with each other. % such that $w_i$ fixes every $e_j$ except those who satisfy $[e_j]^L=E_{r_1+\dots+r_i-1}$, the basis vector of $\t^{L,*}$ coming from the last copy of $A_{l_i-1}$ in $R_L$.
We conclude that $R(\xi)\simeq \Z_2^d$, where $d$ is the number of
type $D$ root systems in $R_0(\xi)$ (including the empty ones, see
\eqref{algemeen}), which is equal to the number of $A$-strips of
length $\kappa_i$ which can be glued to $T_m(\mu)$.\qed

\begin{rem}
Notice that if $l_i=1$ (resp. $l_i=2$), then the corresponding
$L_{l_i}$ are empty (resp. consist of a single root). In these
cases, $w_i \in R(\xi)$ does not correspond to a diagram
automorphism. For example, in the case of minimal principal series
where $\kappa=(1^n)$ we have $R_0(\xi)=R_0$, unless $m=0$ in which
case $R_0(\xi)=D_n$ and hence $R(\xi) \simeq \Z_2$.
\end{rem}

\subsection{Counting irreducible components}

In this section we will show that the induced representation
$\pi(\xi)$ as in the above theorem, decomposes into $2^d$
irreducible components.

Suppose that $R(\xi) \simeq \Z_2^d$. We compute
$H^2(R(\xi),\C^*)$. It is well known that for a direct product of
finite abelian groups $G_1,G_2$, one has \[ H^2(G_1 \times G_2,
\C^*) \simeq H^2(G_1,\C^*) \times H^2(G_2,\C^*) \times (G_1
\otimes G_2).\] Furthermore, for a cyclic group $G$ it is known
that $H^2(G,\C^*)=0$. Since $\Z_2 \otimes \Z_2 \simeq \Z_2$, it is
not hard to see that
\[ H^2(\Z_2^d,\C^*) \simeq \Z_2^{d(d-1)/2}.\]
Let $\eta_\xi$ be the 2-cocycle of Theorem \ref{DOthm}. %We want to show that even though in general $H^2(R(\xi),\C^*) \neq 0$, the number of irreducible summands in $\pi(\xi)$ is equal to $2^p$.
We recall from  \cite[p. 87]{arthur} how to obtain from $\eta_\xi$
a parametrization of the irreducible constituents of $\pi(\xi)$.
One considers a central extension $R'(\xi)$ of $R(\xi)$ such that $\eta_\xi$ splits when pulled back to $R'(\xi)$,
\begin{equation} \label{H}
 1 \to H \to R'(\xi) \to R(\xi) \to 1.
\end{equation}
We may take $H$ to be the cyclic group generated by
$[\eta_\xi]$, the image of $\eta_\xi$ in $H^2(R(\xi),\C^*)$.
Choose a function $\sigma_\xi: R'(\xi) \to \C^*$ which splits
$\eta_\xi$ over $R'(\xi)$, i.e.,
\begin{equation}
\label{sigma}
\eta_\xi(r_1',r_2')=\sigma_\xi(r_1'r_2')\sigma_\xi(r_1')^{-1}\sigma_\xi(r_2')^{-1};
\ r_1',r_2' \in R'(\xi).
\end{equation}

Here we have written, by abuse of notation, $\eta_\xi$ for the
pull-back of $\eta_\xi$ to $R'(\xi)
\times R'(\xi)$. %, i.e., we have $\eta_\xi(r_1',r_2')=\eta_\xi(p(r_1'),p(r_2'))$.
It follows easily that $\sigma_\xi$ satisfies
$\sigma_\xi(zr')=\chi_\xi(z)\sigma_\xi(r')$ for all $z \in H, r'
\in R'(\xi)$. The irreducible characters of $R'(\chi)$ with
central character $\chi_\xi^{-1}$ on $H$ are the ones which index
the irreducible components of $\pi(\xi)$, cf. \cite[p.
87]{arthur}.

Clearly, since $[\eta_\xi] \in H^2(R(\xi),\C^*) \simeq \Z_2^d$, we
either have $H$ is trivial or $H \simeq \Z_2$. If $H$ is trivial,
then we obtain a bijection between the irreducible components of
$\pi(\xi)$ and the irreducible characters of $R(\xi)$, in which
case we are done. Thus, let us consider the second case.

\begin{prop}
Consider a central extension
\begin{equation}
\label{centext}
1 \to H \to G \to \Z_2^n \to 1.
\end{equation}
where $H \simeq \Z_2$. Denote the elements of $H$ by $\pm 1$, and identify them with the corresponding elements in $G$. Then we have, for every irreducible character $\chi \in \hat{G}$, $\chi(-1)=\pm\chi(1)$ and
\begin{equation}
|\{ \chi \in \hat{G} \mid \chi(1)=\chi(-1)\}|=2^n.
\end{equation}
Moreover, all characters in the above set correspond to one-dimensional representations.
\end{prop}

\pf: %Let $H$ be the normal subgroup of $G$ isomorphic to $\Z_2$ which is embedded into the center of $G$ by the central extension.
%Obviously $H$ has only two irreducible characters $\psi_1$ and $\psi_{-1}$ (where $\psi_1$ denotes the trivial character).
By Schur's lemma, the elements of $H$ must act by scalar
multiplication. Therefore the first statement follows.

Since $G$ is a finite group we have
\[ \sum_{\chi \in \hat{G}}\chi(1)\chi(g)=\chi_{\rm reg}(g),\]
where $\chi_{\rm reg}$ is the character of the regular representation of $G$. Applied with $g=1$ and $g=-1$ and summing up, we get
\[ \sum_{\chi \in \hat{G}}\chi(1)(\chi(1)+\chi(-1))=|G|=2^{n+1}.\]
In view of $\chi(1)=\pm \chi(-1)$, we get
\[ \sum_{\chi \in \hat{G}: \chi(-1)=\chi(1)} \chi(1)^2=2^n.\]

By construction $G/H$ is isomorphic to $\Z_2^n$. Thus, by pulling back the $2^n$ irreducible one-dimensional representations of $G/H$ to $H$, we obtain $2^n$ one-dimensional representations of $G$ which are trivial on $H$.

By the above formula, there can be no others and the proposition follows.\qed

\begin{cor}
Consider the central extension \eqref{centext}. Let $N_\pm(G):=|\{ \chi \in \hat{G} \mid \chi(-1)=\pm\chi(1)\}$. Then $N_+(G)=2^n$ and if $N_-(G) \neq 2^n$, then $N_-(G) \leq 2^{n-2}$.
\end{cor}
\pf: The statement on $N_+(G)$ has been proven already. Obviously
$N_-(G)=2^n$ if and only if $G$ is abelian. Let us therefore
suppose that $G$ is not abelian.

Since $[G,G]$ is mapped into $[\Z_2^n,\Z_2^n]=1$, we have $[G,G] \subset {\rm ker}(\pi)={\rm im}(i)$ where $\pi: G \to \Z_2^n$ is the projection map and $i: \Z_2 \to G$ is the inclusion map. Thus $[G,G] \subset \Z_2$ but since $G$ is not abelian, $[G,G] \simeq \Z_2$.

Therefore the exact sequence \eqref{centext} reads $1 \to [G,G] \to G \to G/[G,G] \to 1$ and so $G/[G,G] \simeq \Z_2^n$. The number of one-dimensional representations of a group $G$ being equal to the cardinality of $G/[G,G]$, we see that $G$ has $2^n$ one-dimensional representations and they are exactly those representations for which the restriction to $H$ is trivial.

We get that
\begin{eqnarray*} 2^{n+1} &=&  \sum_{\chi: \chi(1)=\chi(-1)}\chi(1)^2+ \sum_{\chi: \chi(1)=-\chi(-1)}\chi(1)^2 \\ &=& 2^n +  \sum_{\chi: \chi(1)=-\chi(-1)}\chi(1)^2 \\ & \geq & 2^n + 4N_-(G),\end{eqnarray*}
so indeed $N_-(G) \leq 2^{n-2}$ as claimed. \qed

\begin{thm}\label{deco} Consider $\H({\mathcal R},q)$ of type $B_n$, and let $q_1,q_2$ be defined by \eqref{q1q2}. Suppose that $q_2=q_1^m, q_1 \neq 1$.
Let $\Pi_L \subset \Pi_0$ such that $R_L$ is of type $A_\kappa \times B_l$, and consider the
induction datum $\xi=(\Pi_L,\delta,1)$ for $\delta \in \hat{\H}^{ds}_{L,\R}$.
Let $W_Lr_L$ be the central character of $\delta$. Let $\mu \vdash
l$ be such that $W_Lr_L=W_L(r_{\kappa_1}(q_1) \times \dots \times
r_{\kappa_r}(q_1) \times c(\mu;q_1,q_2))$.

Suppose that we can glue $d$ strips of the form $|S(\kappa_i)|$ to
$T_m(\mu)$. Then the induced representation $\pi(\xi)$ decomposes
into a direct sum of $2^d$ inequivalent tempered representations.
\end{thm}

\pf: We will prove the theorem by induction on $d$. If $d=0$ then
$R(\xi)=1$ and thus $\pi(\xi)$ is irreducible. If $d>0$ then
suppose that $R_L$ has type $A_\kappa \times B_l$. Write
$\kappa=(l_1^{r_1}l_2^{r_2}\dots l_s^{r_s})$ (as in
\eqref{hergroep}). We may assume for simplicity of notation that
every $A$-strip $|S(l_i)|$ can be glued to the $m$-tableau
$T_m(\mu)$ (i.e., we assume $s=d$).

Let $\kappa'=(l_2^{r_2}\dots l_p^{r_p})$ and put $n'=l+|\kappa'|$.
Consider the parabolic Hecke algebra (as in \ref{parabolic}) $\H^{L'}$ where $R_{L'}$ has type
$A_{l_1}^{r_1} \times B_{n'}$ and the obvious choice of simple roots $\Pi_{L'}$. Then, by transitivity of induction
we have
\begin{equation}
\label{trans}\pi(\xi)= {\rm Ind}_{\H^L}^\H(\delta)={\rm
Ind}_{\H^{L'}}^\H ({\rm Ind}_{\H^L}^{\H^{L'}}(\delta))={\rm
Ind}_{\H^{L'}}^\H(\pi'(\xi)).
\end{equation}

Clearly the number of irreducible constituents of $\pi(\xi)$ is at
least equal to the number of irreducible constituents of
$\pi'(\xi)$. We want to apply the induction hypothesis to
$\pi'(\xi)$, but we have to circumvent the technical complication
that $\H^{L'}$ is not of type $B$.

We therefore return to the setting of \ref{GHA}. Let $\Hgr^L$
(resp. $\Hgr^{L'}$) be the graded Hecke algebra associated to
$\H^L$ (resp. $\H^{L'}$) as in \ref{GHA}. Thus, $\Hgr^L$ is
associated to the degenerate root datum
$(R_L,\real^*,R_L^\vee,\real,\Pi_L)$ and labels $k_\a$ as in
\eqref{kq}; and $\Hgr^{L'}$ is associated to
$(R_{L'},\real^*,R_{L'}^\vee,\real,\Pi_{L'})$ and labels depending
on those of $\H^{L'}$ as in \eqref{kq}.

By Lusztig's theorems of \cite{lusztig}, we have an equivalence of categories between $\hat{\H}^L_{W_Lr_L}$ (the category of irreducible representations of $\H^L$ with central character $W_Lr_L$) and $\hat{\Hgr}^L_{W_L\gamma_L}$ (the category of irreducible representations of $\Hgr^L$ with central character $W_L\gamma_L$), where ${\rm exp}(\gamma_L)=r_L$. Likewise, we have a categorial equivalence between $\hat{\H}^{L'}_{W_{L'}r_{L}}$ and $\hat{\Hgr}^{L'}_{W_{L'}\gamma_{L}}$ (the central character of any irreducible constituent of $\pi'(\xi)$ is $W_{L'}r_L$).

Let $V$ be the representation space of $\delta$, then $V$ is also a representation space for the $\Hgr^L$-representation corresponding to $\delta$ under the above equivalence. By Lusztig's theorems, the number of irreducible constituents of $\pi'(\xi)$ is equal to the number of irreducible constituents of the corresponding induced representation
\[ {\rm Ind}_{\Hgr^{L}}^{\Hgr^{L'}}(V)=\Hgr^{L'} \otimes_{\Hgr^L} V.\]

Recall that $\Hgr^{L'}=\C[W_0(R_{L'})] \otimes S(\real^*_\C)$ and $\Hgr^L=\C[W_0(R_L)] \otimes S(\real^*_\C)$. Let $L=L_1 \cup L_2$ where $L_1$ consists of those simple roots in $L$ corresponding to the factors $A_{l_1}^{r_1}$ and $L_2$ to the others. Put $\real_{L_1}^*=\{ x \in \real^* \mid \langle x,R_{L_2}^\vee\rangle =0 \}$ and $\real_{L_2}^*=\{ x \in \real^* \mid \langle x, R_{L_1}^\vee \rangle =0 \}$. Then we obtain a decomposition $\Hgr^L=\Hgr^{L_1} \otimes \Hgr^{L_2}$, where $\Hgr^{L_i} (i=1,2)$ is associated to $(R_{L_i},\real^*_{L_i},R_{L_i}^\vee,\real_{L_i},\Pi_{L_i})$ and $k_{L_i}$, the restriction of $k$ to $R_{L_i}$. Therefore, we have a decomposition $V= V_1 \otimes V_2$ where $V_i$ is a representation of $\Hgr^{L_i}$. Likewise, we have a decomposition $\Hgr^{L'}=\Hgr^{L_1} \otimes \Hgr^{n'}$, where $\Hgr^{n'}$ is associated to $(R_0(B_{n'}),\real_{L_2}^*,R_0(B_{n'})^\vee,\real_{L_2},\Pi_{n'})$ and the restriction of $k$ to $R_0(B_{n'})$.

%Clearly, the central character of $V_1$ is $W_0(R_{L_1})(r_{\kappa_1}(k_1) \times \dots \times r_{\kappa_{l_1}}(k_1))$ and the central character of $V_2$ is $W_0(R_{L_2})(r_{\kappa_{l_1+1}}(k_1) \times \dots \times r_{\kappa_r}(k_1) \times c(\mu;k_1,k_2))$.

Therefore, we have
\[ {\rm Ind}_{\Hgr^{L}}^{\Hgr^{L'}}(V)=\Hgr^{L'} \otimes_{\Hgr^{L}}(V)=(\Hgr^{L_1} \otimes \Hgr^{n'})\otimes_{\Hgr^{L_1} \otimes \Hgr^{L_2}}(V_1 \otimes V_2)\]\[ = V_1 \otimes ( \Hgr^{n'} \otimes_{\Hgr^{L_2}} V_2)= V_1 \otimes ({\rm Ind}_{\Hgr^{L_2}}^{\Hgr^{n'}}(V_2)).\]

We now apply again Lusztig's theorems, for the appropriately
defined affine Hecke algebras corresponding to $\Hgr^{L_2}$ and
$\Hgr^{n'}$. Since the latter is an affine Hecke algebra of type
$B$ and the former a parabolic subalgebra, we can apply the
induction hypothesis to the induced representation on the right
hand side. It is easy to see that the central character of the
$\H^{L_2}$-representation affored by $V_2$ is equal to
$r_{\kappa_{r_1+1}}(q_1) \times \dots r_{\kappa_r}(q_1) \times
c(\mu;q_1,q_2)$. Since by construction, we can glue $d-1$ strips
$|S(l_i)| (i=2,\dots,d)$ to $T_m(\mu)$, we may thus assume that
$\pi'(\xi)$ has $2^{d-1}$ irreducible components.

%(...nog verder aan werken...)

%Clearly,  By transitivity of induction and the induction assumption (i.e., that the theorem holds for $p-1$), we know that $\pi(\xi)$ has at least $2^{p-1}$ irreducible components.

Now we consider the number of irreducible components of
$\pi(\xi)$. We have $R(\xi) \simeq \Z_2^d$. If the cocycle
$\eta_\xi$ of Theorem \ref{DOthm} is trivial in $H^2(R(\xi),\C^*)$
then as remarked there, it follows that $\pi(\xi)$ has $2^d$
irreducible components. If not, then (see \cite[p. 86]{arthur}) we
have to consider a central extension $1 \to \Z_2 \to R'(\xi) \to
\Z_2^d \to 1$. Then the irreducible components of $\pi(\xi)$ are
in bijection with the irreducible representations of $G$ with a
fixed restriction to $\Z_2$, as explained above. We denote this
character of $\Z_2$ by $\chi$.

Recall that $N_+(R'(\xi))=2^d$ and if $N_-(R'(\xi)) \neq 2^d$ then
$N_-(R'(\xi))\leq 2^{d-2}$. Since we know that $\pi(\xi)$ has at
least $2^{d-1}$ irreducible components, it follows in this case as
well that $\pi(\xi)$ is a direct sum of $2^d$ irreducible
components, which are indexed by the irreducible representations
of $R(\xi)$.

These components are tempered by \cite[Thm 4.23]{O1}. The fact
that they are mutually inequivalent follows since the multiplicity
of the constituent of $\pi(\xi)$ which is indexed by $\rho \in
\hat{R}(\xi)$ is equal to ${\rm dim}(\rho)=1$, see again \cite[p.
87]{arthur}. \qed %Indeed, as in \cite[p. 87]{arthur}, one obtains a representation of $R'(\xi) \times \H$ on $\pi(\xi)$. This representation decomposes into a direct sum
%\[ \oplus_{\rho} m_\rho (\rho^\vee \otimes \pi_\rho).\] In this sum, $\rho$ ranges over the irreducible characters of $R'(\xi)$. The multiplicity $m_\rho$ is zero if $\rho$ does not restrict to $\chi$ on $\Z_2$, and one if it does. As $\rho$ ranges over the characters in the last set, $\pi_\rho$ ranges over the irreducible constituents of $\pi(\xi)$, counting all exactly once. Thus, the multiplicity of $\pi_\rho$ is equal to the dimension of $\rho$.  The statement follows since we have seen that the irreducible representations of $R'(\xi)$ which index the constituents of $\pi(\xi)$ are all one-dimensional.\qed

%It is well-known (since \cite{schur1911}) that $H^2(\Z/2\Z,\C^*)=0$. Therefore, the cocycle $\eta_\xi$ always splits and by the remarks in paragraph \ref{decomp}, the irreducible constituents of $\pi(\xi)$ are in bijection with the irreducible characters of $R(\xi)$. Since $R(\xi) \simeq (\Z/2\Z)^x$, the number of irreducible characters of $R(\xi)$ is $2^x$.\qed
%\vspace{5cm}

In turn, it follows that $\eta_\xi$ must split:

%In general, the map $R(\xi) \to {\rm End}_\H(\pi(\xi)): r \mapsto E(r,\xi)$ needs not be a homomorphism: one only has $E(r_1,\xi)E(r_2,\xi)=\eta_\xi(r_1,r_2) E(r_1r_2,\xi)$ for a certain 2-cocycle $\eta_\xi: R(\xi) \times R(\xi) \to \C^*$.

\begin{prop} The 2-cocycle $\eta_\xi$ has trivial image in
$H^2(R(\xi),\C^*)$.
\end{prop}

\pf: Suppose this is not true, then the group $H$ in \eqref{H} is
isomorphic to $\Z_2$.

Suppose that $R'(\xi)$ is not abelian. Then, as we
have seen in the proof of Theorem \ref{deco}, $\chi_\xi$ must be
the trivial character. But this means that $\sigma_\xi$, the
function of \eqref{sigma} which splits $\eta_\xi$ on $R'(\xi)$,
satisfies $\sigma_\xi(zr')=\sigma_\xi(r')$ for all $z \in H, r'
\in R'(\xi)$. %Define a function, also denoted $\sigma_\xi$, on $R(\xi)$ by $\sigma_\xi(r):=\sigma_\xi(p(r))$. Since $p(r_1r_2)=zp(r_1)p(r_2)$ for some $z \in H$, we obtain that
So $\sigma_\xi$ descends to $R(\xi)$ and splits $\eta_\xi$. In
other words, the image of $\eta_\xi$ in $H^2(R(\xi),\C^*)$ is
trivial, which is a contradiction.

On the other hand, suppose that $R'(\xi)$ is abelian. The equivalence classes of extensions of the form \eqref{H} are in a natural bijection with the elements of $H^2(R(\xi),\Z_2)$. %which is isomorphic to $\Z_2^{p(p+1)/2}$ if $R(\xi) \simeq \Z_2^p$. The abelian extensions of the form \eqref{H} correspond under this bijection to ${\rm Ext}(R(\xi),\Z_2)$.
Therefore, under the map $H^2(R(\xi),\Z_2) \to
H^2(R(\xi),\C^*)$, the cohomology class of the central extension
\eqref{H} is mapped to $[\eta_\xi] \in {\rm Ext}(R(\xi),\C^*)$.
 %For a finite cyclic group $G$ of order $n$ and an abelian group $A$, it is known that ${\rm Ext}(G,A) \simeq A/A^n$ (see e.g. \cite[p. 52]{karp85}), which implies that in particular, ${\rm Ext}(\Z_2,\Z_2)=0$.
%Notice that ${\rm Ext}(\Z_2,\Z_2)=\Z_2$ since $H^2(\Z_2,\Z_2)=\Z_2$ and there are only abelian extensions of $\Z_2$ by $\Z_2$. Now we use that for abelian groups $G_1,G_2,A$ one has ${\rm Ext}(G_1 \times G_2, A) \simeq {\rm Ext}(G_1,A) \times {\rm Ext}(G_2,A)$ (see e.g. \cite[p. 63]{karp85}) to conclude that ${\rm Ext}(R(\xi),\Z_2) \simeq R(\xi)$. %It remains to identify the image of ${\rm Ext}(R(\xi),\Z_2)$ in $H^2(R(\xi),\C^*)$. %If we choose 2-cocycles representing the elements in ${\rm Ext}(R(\xi),\Z_2)$, then these are precisely the symmetric ones: those 2-cocycles $f: R(\xi) \times R(\xi) \to \Z_2$ for which $f(r_1,r_2)=f(r_2,r_1)$ for all $r_1,r_2 \in R(\xi)$.
%However, ${\rm Ext}(\Z_2,\C^*)=0$ (since $H^2(\Z_2,\C^*)=0$),
%hence also ${\rm Ext}(R(\xi),\C^*)=0$ since for abelian groups
%$G_1,G_2,A$ one has ${\rm Ext}(G_1 \times G_2,A) \simeq {\rm
%Ext}(G_1,A) \times {\rm Ext}(G_2,A)$ (see e.g. \cite[p.
%63]{karp85}). Thus, $[\eta_\xi]$ is trivial, %$[\tilde{\eta}_\xi] \mapsto 0$ under the map $H^2(R(\xi),\Z_2) \to H^2(R(\xi),\C^*)$. This means that $[\eta_\xi]=0$,
%which contradicts our assumption. The Proposition follows. \qed
Since $\C^*$ is a divisible abelian group, it is injective. Hence ${\rm Ext}(R(\xi),\C^*)=1$ and in particular, $[\eta_\xi]=1$.\qed

\subsubsection{A combinatorial remark}
By Theorem \ref{deco}, the number of irreducible constituents of
an induced discrete series representation with real central
character is as predicted by Conjecture \ref{con}. Moreover, the
parametrization of the constituents of $\pi(\xi)$ by the
irreducible characters of $R(\xi)$ can be transferred to a
parametrization in terms of Young tableaux. Indeed, consider
$\pi(\xi)$ where $\xi=(\Pi_L,\delta,1) \in \Xi$ has real central
character $W_Lr_L$ with $r_L$ as in \eqref{rL}. Then the number of
irreducible components in $\pi(\xi)$ is $2^d$ where $d$ is the
number of strips of length $\kappa_i$ which can be glued to
$T_m(\mu)$. For every $J\subset \{1,\dots,d\}$, we obtain a
partition $\mu_J$, consisting of $T_m(\mu)$ with the strips
corresponding to $J$ glued to it. We obtain $2^d$ Young tableau
associated to $\pi(\xi)$. Let $\chi_J$ be the character of
$R(\xi)$ which coincides with the non-trivial character of $\Z_2$
on the factors which correspond to $J$ and which coincides with
the trivial representation on the other factors. Then we have a
natural bijection $\chi_J \leftrightarrow \mu_J$. Since $\chi_J$
corresponds to an irreducible component of $\pi(\xi)$, we can also
parametrize this component by $\mu_J$. In this way we obtain a
parametrization of the irreducible constituents of $\pi(\xi)$ by
the $m$-tableaux $T_m(\mu_J)$.

%These tableaux arise naturally in pairs $\mu_J$ and $\mu_{J^c}$,
%where $J^c$ denotes the complement of $J$ in $\{1,\dots,p\}$. In
%terms of Springer correspondence, we have the following
%incarnation of these pairs. Let $\kappa_J$ denote the partition
%with parts $l_j, j\in J$. Let $L_J=r_{L_J}T^{L_J}$ be the residual
%coset for which $R_L$ has type $A_{\kappa-\kappa_J} \times
%B_{|\l_J|}$ and whose center, the (generically) residual point
%$r_{L_J}$ coincides on the $B_{|\mu_J|}$-factor with
%$c(\mu_J;q_1,q_2)$. Choosing $\e>0$ sufficiently small, it has
%been shown in \cite{proefschrift} that
%\[ \Sigma_{m\pm\e}(W_0r_{L_J})=\Sigma_{m\mp\e}(W_0r_{L_{J^c}}).\]
%In other words, the Springer correspondents of the central
%characters $W_0r_{L_J}$ and $W_0r_{L_{J^c}}$ are exchanged when we
%compare the Springer correspondence for $q_2=q_1^{m+\e}$ with the
%one for $q_2=q_1^{m-\e}$. In terms of the irreducible characters
%of $R(\xi)$, one may view this as an exchange between the
%characters $\chi_J$ and $\chi_{J^c}$.

\subsubsection{The Iwahori-Hecke cases}
Let all $q_{\ca}=q$, the cardinality of the residue field of the $p$-adic field $F$. Then $\H \simeq {\rm End}_G({\rm Ind}_\I^G(1))$ for the group $G=G(F)$ with root datum $(R_0^\vee,Y,R_0,X,\Pi_0^\vee)$. It is known from the work of Kazhdan and Lusztig, \cite{KL}, that the restrictions to $\H_0$ of the modules in $\Hrcc$ specialize, for $q_{\ca}^{1/2} \to 1$, into the Springer modules for the Weyl group of the complex group $\G$ with root datum $(R_0,X,R_0^\vee,Y,\Pi_0)$. For example, if $R_0$ is of type $A_n$ (resp. $B_n$, resp. $C_n$) then $\G=SL_{n+1}(\C)$ (resp. ${\rm Spin}_{2n+1}(\C)$, resp $Sp_{2n}(\C)$). It is also known that we obtain a bijection between the central characters of the modules in $\Hrcc$ and the unipotent conjugacy classes of $\G$. Since there is a canonical bijection between the unipotent conjugacy classes of $\G$ and those of $\G/Z$ (where $Z$ denotes its center) we replace $\G$ with the corresponding group of adjoint type. We denote the bijection between the central characters of $\Hrcc$ and the unipotent classes of $\G$ by $W_0r \mapsto u_{W_0r}$ where $u_{W_0r}$ is a representative of the unipotent class associated to $W_0r$. Let, for $u \in \G$, $A(u)$ be the group of connected components of the centralizer $C_\G(u)$.

In our situation where $R_0$ is of type $B_n$ we obtain two of these cases, cf. remark \ref{BC}. If $m=1$, then the modules of $\Hrcc$ specialize into the Springer modules of $\G=SO_{2n+1}(\C)$. If $m=1/2$, then the modules of $\Hrcc$ are in bijection with those of an affine Hecke algebra with root system of type $C_n$ and equal labels. In this case, the modules of $\Hrcc$, restricted to $\H_0$, specialize into the Springer modules of the adjoint group ${\rm Spin}_{2n}(\C)$.

In this context, the $R$-group admits the following characterization:
\begin{prop}\label{Au} Let $m=1$ or $m=1/2$. Let $\xi=(\Pi_L,\delta,1) \in \Xi$.
Let $W_Lr_L$ be the central character of $\delta \in \hat{\H}_{L,\R}^{ds}$. Then
\[ R(\xi) \simeq A(u_{W_0r_L})/A(u_{W_Lr_L}).\]
\end{prop}

\pf: First we recall the relation between $A(u_{W_0r_L})$ and the $m$-symbols of  the Springer correspondents $\Sigma_m(W_0r_L)$, for a central character $W_0r_L$ of a representation in $\Hrcc$. Let $S$ be the set of entries which occur exactly once in the symbol $\Sigma_m(W_0r_L)$. Lusztig has defined in \cite{luscells} (for an overview of these results, see \cite[p. 419]{carter}) an interval in $S$ to be a subset of the form $S=(i,i+1,\dots,j)$, such that $0 \leq i \leq j$, $i-1 \notin S, j+1 \notin S$. If $m=1/2$, then we also require $i \geq 1$. Consider the group $A'(u)=\Z_2 \times \dots \times \Z_2$, one copy of $\Z_2$ for every interval in $S$. Let $a_i$ generate the $i$-th copy of $\Z_2$. Then, if $m=1$,  $A(u) \simeq \{ \sum_i n_ia_i \mid \sum_i n_i {\rm\ even}\} \subset A'(u)$. If $m=1/2$ then we need to consider $\l(u)$ such that $u \in {\rm Spin}_{2n}(\C)$ has elementary divisors partition $\l(u)$. Let $\l(u)=(1^{r_1}2^{r_2}\dots 2n^{r_{2n}})$. Then the number of intervals is equal to the number of $r_i$ such that $i$ is even and $r_i>0$. We associate to every interval a generator $a_i$ of $A'(u)$ and a $r_i$, using the natural ordering. Then $A(u)=A'(u)/ \sum_{i\ {\rm even}, r_i \ {\rm odd}}a_i$. %Conceptually, in both cases, the intervals are the subsets of $S$ whose lowest entry  can both appear in the lower and in the upper row of two similar $m$-symbols. (niet helemaal waar, bv als er maar 1 entry is)

Now let $W_Lr_L$ be as in the statement, and write it in the form \eqref{rL}.  Observe that $A(u_{W_Lr_L})=A(u_{W_0(B_l)c(\mu;q_1,q_2)})$, since $A(u)=1$ for all $u$ in groups of type $A$. Suppose that $d$ strips $|S(\kappa_i)|$ can be glued to $T_m(\mu)$. Then $R(\xi) \simeq \Z_2^d$ and we have to show that the number of intervals in the $m$-symbols of $\Sigma_m(W_0r_L)$ is exactly $d$ greater than the number of intervals in the $m$-symbols of $\Sigma_m(W_0(B_l)c(\mu;q_1,q_2))$, where for $m=1/2$ we have to also check that $\l(u_{W_0r_L})$ has even parts with odd multiplicity if and only if $\l(u_{W_0(B_l)c(\mu;q_1,q_2)})$ has even parts with odd multiplicity. But this is obvious since $\l(u_{W_0r_L})$ is obtained from $\l(u_{W_0(B_l)c(\mu;q_1,q_2)})$ by adding the parts $(\kappa_1,\kappa_1,\kappa_2,\kappa_2,\dots,\kappa_r,\kappa_r)$ (this follows from the Bala--Carter classification, see e.g. \cite{jantzen} for a nice presentation).

Let $W_0r$ be a residual point. Then we have seen in \cite{S} that the entries of $\Sigma_m(W_0r)$ all have a difference of at least two, hence every entry forms its own interval, except for the entry zero and $m=1/2$. %(Note that if $W_0r=W_0c(\mu;q_1,q_2)$ where $T_m(\mu)$ has $2r$ or $2r+1$ blocks, then $A(u) \simeq \Z_2^{2r}$, and there are ${2r+1 \choose r}$ symbols in the similarity class of $\Sigma_m(W_0r_L)$, so $A(u)$ has this number of geometric characters).

 Then we need to consider the symbols of the Springer correspondents $\Sigma_m(W_0r_L)$. Since truncated induction is transitive (\cite[Cor. 4.35]{S}) , we define $(\alpha^{(0)},\beta^{(0)})=\S_m(\mu)$ and we choose $(\a^{(i)},\beta^{(i)}) \in {\rm tr}_m-{\rm Ind}((\kappa_i) \otimes (\a^{(i-1)},\beta^{(i-1)}))$ for $i=1,\dots,r$. By \cite[Prop. 4.37]{S}, we have $(\a^{(i)},\beta^{(i)})=(\a^{(i-1)},\beta^{(i-1)}) \cup (a_i,b_i)$ where $a_i+b_i=\kappa_i$. The similarity class of $\Sigma_m(\a^{(i)},\beta^{(i)})$ is independent of this choice. Let $S^{(i)}$ denote the number of intervals in the $m$-symbol of $(\a^{(i)},\beta^{(i)})$.

We recall the following from \cite[Props. 4.41-4.43]{S} and their proofs. If $|S(\kappa_i)|$ can not be glued to $T_m(\mu)$, then $(\a^{(i)},\beta^{(i)})$ is uniquely determined by the choice of $(\a^{(i-1)},\beta^{(i-1)})$, and the entries of $a_i,b_i$ in the $m$-symbol $\Lambda^m(\a^{(i)},\beta^{(i)})$ are either equal, or have a difference of one. Suppose that they are equal. Since no entry can occur more than twice, it follows that $|S^{(i)}|=|S^{(i-1)}|$. If the entries of $a_i,b_i$ are not equal, then they have a difference of one, and moreover form an interval together with exactly one of the other already existing intervals. Thus, in this case $|S^{(i)}|=|S^{(i-1)}|$ as well. On the other hand, if one can glue a strip $|S(\kappa_i)|$ to $T_m(\mu)$, then the entries of $a_i,b_i$ in the symbol of $(\a^{(i)},\beta^{(i)})$ have a difference of one and do not form an interval with the other entries. Thus, in this case we get $|S^{(i)}|=|S^{(i-1)}|+1$.

In total, the number of intervals in the $m$-symbol of any element of $\Sigma_m(W_0r_L)$ is therefore equal to $|S^{(r)}|=|S^{(0)}|+d$, where $d$ is the number of strips $|S(\kappa_i)|$ which can be glued to $T_m(\mu)$ and $|S^{(0)}|$ is the number of intervals in the $m$-symbol of $\Sigma_m(W_0(B_l)c(\mu;q_1,q_2))$. Therefore, the result follows.\qed

\begin{ex}  We continue the example of Figure \ref{example}, even though $m>1$ there. However, even without the interpretation in terms of component groups, the statement on the number of intervals needed in the proof remains true, if we define for $m \in\Z$ an interval to be a maximal set of consecutive entries in the set of entries which occur precisely once. In the example, we have $\S_3(\mu)=(234,2)$ which has 3-symbol
\[ \left( \begin{array}{ccccccc}0&&4&&7&&14\\ &2&&&&& \end{array}\right).\]
As remarked in the proof of \ref{Au}, this symbol has 5 intervals, one for every entry.
We now perform the truncated induction ${\rm tr}_3-{\rm Ind}((3) \otimes (4) \otimes (7) \otimes (11) \otimes (234,2))$. One computes, using the methods of \cite{S}, that the $m$-symbols of the 2-partitions which occur in this induction are all in the similarity class of
\[ \left( \begin{array}{ccccccccccccccc} 0&&2&&6&&8&&11&&13&&16&&18 \\ &1&&4&&6&&10&&15&&&& \end{array} \right)\]
This symbol has 7 intervals, namely $(0,1,2),(4),(8),(10,11),(13),(15,16),(18)$. Indeed, this is two more than the number of intervals in the symbol of $\S_3(\mu)$, and one can glue two strips to $T_3(\mu)$; those of length 7 and 11. These correspond to the intervals $(10,11)$ and $(15,16)$. The strips of length 3 and 4 correspond to the interval $(0,1,2)$ (the extension of an already existing one) and to the numbers $6,6$ which do not appear in any interval.

\end{ex}

%Let $l \geq 0$ and consider a 2-partition $(\alpha,\beta)$ of $l$. We have seen in \cite{S} that if we consider the truncated induction \eqref{trind} of the form ${\rm tr}_m-{\rm Ind}_{S_{\kappa_i} \times W_0(B_l)}^{W_0(B_{\kappa_i+l})}((\kappa_i) \otimes (\alpha,\beta))$

\medskip
\noindent{\bf Acknowledgments.} The author wishes to thank Eric
Opdam for useful discussions and comments on earlier versions of
this paper.

\bibliographystyle{amsplain}
\bibliography{biblio}
%\thebibliography

\end{document}